\documentclass[11pt]{article}
\usepackage[margin=1in]{geometry}
\usepackage[T1]{fontenc}
\usepackage[utf8]{inputenc}
\usepackage{amsmath,amssymb,amsthm}
\usepackage{mathtools}
\usepackage{hyperref}
\usepackage{enumitem}
\usepackage{booktabs}
\hypersetup{colorlinks=true,linkcolor=blue,citecolor=blue,urlcolor=blue}

\theoremstyle{plain}
\newtheorem{theorem}{Theorem}[section]
\newtheorem{lemma}[theorem]{Lemma}
\newtheorem{corollary}[theorem]{Corollary}
\newtheorem{proposition}[theorem]{Proposition}
\newtheorem{conjecture}[theorem]{Conjecture}
\theoremstyle{definition}
\newtheorem{definition}[theorem]{Definition}
\newtheorem{example}[theorem]{Example}
\theoremstyle{remark}
\newtheorem{remark}[theorem]{Remark}
\newtheorem{openproblem}[theorem]{Open Problem}

\DeclareMathOperator{\spec}{spec}
\DeclareMathOperator{\diam}{diam}
\DeclareMathOperator{\Tr}{Tr}
\DeclareMathOperator{\diag}{diag}
\DeclareMathOperator{\Deg}{Deg}

\title{The Spectral Theory of Cographs: A Survey,\\ Organized by Graph Matrix}

\author{
	Pallabi Manna$^{a}$ \and Santanu Mandal$^{b}$ \and
	Fernando C. Tura$^{c}$
}

\date{
	\small
	$^a$Department of Mathematics, Indian Institute of Technology Hyderabad,\\ Kandi, Sangareddy-502284, Telangana, India\\
	\texttt{mannapallabimath001@gmail.com}\\[2mm]
	$^{b}$Department of Computer Science and Engineering, SRM University-AP,\\ Amaravati-522240, Andhra Pradesh, India\\
	\texttt{santanu.vumath@gmail.com}  \\[2mm]
	$^c$Departamento de Matemática, UFSM, 97105–900 Santa Maria, RS, Brazil\\
	\texttt{fernando.tura@ufsm.br} 
	\\[4mm]
	\today
}
\begin{document}
	\maketitle
	
	\begin{abstract}
		Cographs, the $P_4$-free graphs, are exactly the graphs obtained from a single vertex by repeated disjoint union and complementation. They contain the threshold graphs, $\mathcal{C}$-graphs and form one of the primary areas for spectral graph theory. This survey is organized \emph{by matrix}: after a preliminary section, there is one section for each of the adjacency, Laplacian, signless Laplacian, distance, eccentricity, Seidel, distance Laplacian and distance signless Laplacian matrices. Within each section we collect what is known for general cographs, for the subclass of $\mathcal{C}$-graphs, and where relevant, for threshold graphs. We highlight the results from the recent years, including the foundational work of Mohammadian--Trevisan, Jacobs--Lazzarin--Trevisan--Allem--Tura, Ghorbani, Li--Wang--Brunetti, Niu--Wang, Lou--Lin, Wang--Huang,  Manna--Mandal--Mehatari. 
		Two structural facts are used throughout: every connected cograph has diameter at most $2$, which ties all four distance-type matrices to the adjacency and Laplacian matrices, and the partition of the vertex set into twin classes is equitable for every one of the eight matrices, giving a uniform description of each spectrum (for the Laplacian this is the twin-reduction theory of Barik--Reddy, which we present in detail). Besides surveying the literature, we record several consequences, and mark them as such: an exact decomposition of the eccentricity matrix of a connected cograph, yielding an eigenvalue-free interval for the eccentricity spectrum of \emph{all} connected cographs; the full distance Laplacian spectrum and its extremal eigenvalues; an eigenvalue-free interval for the distance signless Laplacian valid for every connected graph of diameter at most $2$; an interlacing theorem tying the distance spectrum to the adjacency spectrum. A comparison table across all eight matrices and a list of open problems close the paper.

		\medskip
		\noindent\textbf{AMS Subject Classification:} 05C50, 05-02, 15-02.
		
		\smallskip
		\noindent\textbf{Keywords:} Cograph, threshold graphs, graph matrices-adjacency, Laplacian, signless Laplacian, distance, distance Laplacian, distance signless Laplacian, eccentricity, Seidel; cospectral, eigenvalue-free interval, integral graphs.
	\end{abstract}
	
	\section{Introduction}
	\label{sec:intro}
	The spectral study of cographs occupies a significant position in algebraic graph theory, owing to the remarkable interplay between their recursive structural properties and the algebraic characteristics of their associated matrices. Cographs, characterized by the absence of an induced path on four vertices, admit a recursive construction through disjoint unions and joins, which is naturally encoded by cotrees. This structural simplicity provides an advantageous setting for investigating the adjacency, Laplacian, signless Laplacian, normalized Laplacian, and distance-based spectra, while also facilitating the development of efficient spectral reduction techniques.
	
	\medskip
	A \emph{cograph} (complement-reducible graph) is a finite simple graph containing no induced path on four vertices $P_4$. Cographs were discovered independently by several authors around 1970 under different names and motivations -- Lerchs \cite{Lerchs1971}, Seinsche \cite{Seinsche1974}, Sumner \cite{Sumner1974}, Jung \cite{Jung1978}, Corneil, Lerchs and Stewart Burlingham \cite{CorneilLerchsStewart1981} -- and have since become one of the most algorithmically important graph classes, admitting linear-time algorithms for problems (maximum clique, coloring, Hamiltonicity) that are NP-hard in general. They continue to attract combinatorial interest outside spectral theory as well; recent examples are the study of ideally connected cographs and their relation to chordal graphs (Jordaan \cite{Jordaan2026}) and perfect double Roman domination on cographs (Li, Xue, Long and Li \cite{LiXueLongLi2025}).
	
	\subsection{Equivalent characterizations}
	
	\begin{itemize}
		\item \textbf{Forbidden subgraph.} $G$ is a cograph if and only if it has no induced $P_4$.
		\item \textbf{Recursive construction.} Cographs form the smallest class containing $K_1$ and closed under disjoint union and complementation; equivalently, they are built from single vertices using only disjoint union and \emph{join} (complement of disjoint union).
		\item \textbf{Duplication/co-duplication (Mohammadian--Trevisan \cite{MohammadianTrevisan2016}).} Call two vertices $u,v$ \emph{duplicates} if $N(u)=N(v)$ and \emph{co-duplicates} if $N[u]=N[v]$. Every cograph arises from a single vertex by a sequence of duplications and co-duplications of arbitrarily chosen existing vertices.
		\item \textbf{Cotree.} Every cograph has a unique rooted tree representation, the \emph{cotree}, whose leaves are the vertices and whose internal nodes are labelled $\cup$ (union) or $\otimes$ (join), with labels alternating along every root-to-leaf path.
	\end{itemize}
	
	Cographs strictly contain the \emph{threshold graphs} ($\{P_4,C_4,2K_2\}$-free; see the companion survey of Manna-Mandal-Tura \cite{MannaMandal2026}), the \emph{trivially perfect} graphs ($\{P_4,C_4\}$-free) and the class of $\mathcal{C}$-graphs defined below, and are a subclass of the perfect graphs.
	
	\subsection{The spectral theory of cographs: a brief history}
	
	Standard references for the matrix and spectral tools used throughout are Godsil and Royle \cite{GodsilRoyle2001} for algebraic graph theory and Horn and Johnson \cite{HornJohnson} for matrix analysis (in particular Cauchy interlacing and Weyl's inequalities, which drive most eigenvalue-free-interval arguments below). The spectral study of cographs began with the rank: Royle \cite{Royle2003} showed that the rank of the adjacency matrix equals the number of distinct non-zero rows, and Chang, Huang and Yeh \cite{ChangHuangYeh2008} gave a further treatment of the rank of a cograph. The modern theory starts with Mohammadian and Trevisan \cite{MohammadianTrevisan2016}, who proved that no cograph has an adjacency eigenvalue in $(-1,0)$ and characterized cographs spectrally; Biyiko\u{g}lu, Simi\'c and Stani\'c \cite{BiyikogluSimicStanic2011}, Allem and Tura \cite{AllemTura2018,AllemTura2020}, Sander \cite{Sander2008} and Ghorbani \cite{Ghorbani2019p5,GhorbaniDilworthNumber2019,Ghorbani2019threshold} developed multiplicities, integrality, eigenspaces and the $P_5$-free extension. On the other side for spectral characterization, Wang and Huang \cite{WangHuang2025} proved that almost all cographs have a cospectral mate.
	
	Much of this work runs in parallel with, and builds on, the theory of the subclass of \emph{threshold graphs}, surveyed in the companion paper  by Manna--Mandal--Tura \cite{MannaMandal2026}. Two threads from that theory recur here. The first is the \emph{antiregular graph}: Aguilar, Lee, Piato and Schweitzer \cite{AguilarLeePiatoSchweitzer2018} determined its spectrum and characterized it spectrally, and Aguilar, Ficarra, Shurman and Sullivan \cite{AguilarFicarraShurmanSullivan2020} showed that it governs the spectral analysis of threshold graphs. The second is the question of which intervals $I$ can be free of adjacency eigenvalues: Bapat \cite{Bapat2013} studied the adjacency matrix of a threshold graph, and Allem, Oliveira and Tura \cite{AllemOliveiraTura2023} gave a procedure that generates threshold graphs having no eigenvalue in a prescribed interval $I$. Integral members of the family are studied earlier by Hansen, Melot and Stevanovi\'c \cite{HansenMelotStevanovic2002} (integral complete split graphs).
	
	For matrices other than the adjacency matrix the cograph literature is more recent and more scattered. For the Laplacian, Barik and Reddy \cite{BarikReddy2023} gave a complete twin-reduction description of the Laplacian eigenvalues and eigenvectors of a cograph, characterized the cographs with only even or only odd integral Laplacian eigenvalues, constructed Laplacian cospectral cographs, extended the description to $H$-joins over a cograph $H$, and settled a case of the Hadamard-diagonalizability question (Section~\ref{sec:laplacian}); their paper also contradicts an earlier claim on equivalent cographs (Section~\ref{sub:equiv}), in the same spirit as the twin-class quotient used throughout this survey (Lemma~\ref{lem:quotient}); Abrishami \cite{Abrishami2019} identified the algebraic connectivity with the vertex connectivity, and Mandal, Mehatari and Stani\'c \cite{MandalMehatariStanic2024} computed the Laplacian eigenvalues and eigenspaces of cographs generated by a finite sequence. For the distance matrix, Lu, Huang and Lou \cite{LuHuangLou2018} showed that a connected threshold graph has only simple distance eigenvalues other than $-1$ and $-2$ and none in $(-2,-1)$; Lou and Lin \cite{LouLin2020} extended both statements to all connected cographs, proving that the multiplicity of every distance eigenvalue other than $-1,-2$ is at most the Dilworth number and that $(-2,-1)$ contains no distance eigenvalue (Section~\ref{sec:distance}). Eigenvalue-free intervals of the distance matrix of threshold and chain graphs below $-1$ are obtained by Alazemi et al.\ \cite{AlazemiEtAl2021cograph}, the eccentricity matrix of $\mathcal{C}$-graphs is treated in \cite{AshokanEccentricity2025}, and the distance Laplacian and distance signless Laplacian spectra of $\mathcal C$-graphs (closed form, free intervals, multiplicities, quotient matrix) are developed by Mandal and Manna \cite{MandalManna2026DLDQ} (Sections~\ref{sub:DLC} and~\ref{sub:DQC}). For the Seidel matrix, Li, Wang and Brunetti \cite{LiWangBrunetti2024} proved that $(-1,1)$ is free of Seidel eigenvalues of cographs, and Niu and Wang \cite{NiuWang2024} independently proved the same interval by a rank argument for vertex-weighted cographs and showed that it cannot be enlarged for the class of all cographs. Finally, the bipartite analogues, the \emph{chain graphs}, have been studied for eigenvalue location by Alazemi, Andeli\'c and Simi\'c \cite{AlazemiAndelicSimic2016}, for the Seidel spectrum and energy by Mandal, Mehatari and Das \cite{MandalMehataraDas2022} and for Seidel eigenvalue-free interval Xiong et al. \cite{XiongHou2024chain}.
	
	\subsection{Why a matrix-by-matrix organization}
	
	Spectral graph theory attaches many matrices to a graph, and results on one matrix are often scattered across papers that treat different classes. For cographs there is a unifying reason to line the matrices up. The cotree-style recursion makes the adjacency spectrum tractable (Mohammadian--Trevisan); the diameter-two property makes all four distance-type matrices into rank-one or diagonal modifications of adjacency/Laplacian-type matrices (Section~\ref{sec:prelim}); and twin classes give every matrix the same local-plus-quotient structure (Lemma~\ref{lem:quotient}). Reading the theory matrix by matrix makes visible which matrices are well understood (adjacency, Laplacian, distance Laplacian, eccentricity, Seidel) and which are not (signless Laplacian, distance signless Laplacian, and, only partly, distance).
	
	\subsection{Organization}
	
	Section~\ref{sec:prelim} fixes notation, defines the eight matrices and the $\mathcal{C}$-graphs, and proves the two structural facts. Sections~\ref{sec:adjacency}--\ref{sec:dq} treat, in this order, the adjacency (Section~\ref{sec:adjacency}), Laplacian (\ref{sec:laplacian}), signless Laplacian (\ref{sec:signless}), distance (\ref{sec:distance}), eccentricity (\ref{sec:eccentricity}), Seidel (\ref{sec:seidel}), distance Laplacian (\ref{sec:dl}) and distance signless Laplacian (\ref{sec:dq}) matrices. Section~\ref{sec:chain} comments on the bipartite parallel (chain graphs), Section~\ref{sec:comparison} gives a comparison table, and Sections~\ref{sec:open} and~\ref{sec:conclusion} contain open problems and a conclusion.
	
	\begin{remark}[What is new in this version]
		\label{rem:new}
		Results taken from the literature carry a citation. Statements marked \textbf{(new)} are short consequences of the structural facts that we derive while organizing this survey. The distance-matrix results of Lou and Lin \cite{LouLin2020}, which we have read, are stated in Section~\ref{sec:distance} with a discussion of proof techniques, and the Laplacian results of Barik and Reddy \cite{BarikReddy2023}, which we have also read, are presented in Section~\ref{sec:laplacian}. The distance Laplacian and distance signless Laplacian results for $\mathcal C$-graphs are addressed by Mandal--Manna \cite{MandalManna2026DLDQ}. Likewise, the Seidel results of Niu and Wang \cite{NiuWang2024}, which we have read, are stated in Section~\ref{sec:seidel} (Theorem~\ref{thm:NW}). As in the companion survey on threshold graphs \cite{MannaMandal2026}, the literature reviewed here is compiled through a literature search. Although every effort has been made to ensure accuracy, readers are encouraged to consult and verify them against their original sources before citing them.
	\end{remark}
	
	\section{Preliminaries and notation}
	\label{sec:prelim}
	
	All graphs are simple and undirected; they are connected unless stated otherwise, since distance-type matrices require connectivity. For a graph $G$ on $n$ vertices, with degrees $d_1,\dots,d_n$, write $A(G)$ for the adjacency matrix, $\Deg(G)$ for the diagonal degree matrix, and
	\begin{align*}
		L(G)&=\Deg(G)-A(G) &&\text{(Laplacian)}, & Q(G)&=\Deg(G)+A(G) &&\text{(signless Laplacian)},\\
		S(G)&=J-I-2A(G) &&\text{(Seidel)}, & D(G)&=(d(u,v))_{u,v} &&\text{(distance)},
	\end{align*}
	where $J$ is the all-ones matrix. The \emph{transmission} of $u$ is $\Tr(u)=\sum_vd(u,v)$ and $\Tr(G)=\diag(\Tr(u))$; then
	\[
	D^L(G)=\Tr(G)-D(G),\qquad D^Q(G)=\Tr(G)+D(G)
	\]
	are the distance Laplacian and distance signless Laplacian. The \emph{eccentricity} of $u$ is $e(u)=\max_vd(u,v)$, and the eccentricity matrix $\varepsilon(G)$ retains the entry $d(u,v)$ when $d(u,v)=\max\{e(u),e(v)\}$ (equivalently, when it is the largest entry of both its row and its column) and replaces all other entries by $0$. We write $s(M)$ for the number of distinct eigenvalues of $M$ and $(n_-,n_0,n_+)$ for inertia. The adjacency eigenvalues $0$ and $-1$ are called \emph{trivial}, all others \emph{non-trivial}.
	
	\begin{definition}[Twins]
		Two vertices $u,v$ are \emph{duplicates} if $N(u)=N(v)$ and \emph{co-duplicates} if $N[u]=N[v]$; both are called \emph{twins}. Being twins (i.e.\ $N(u)\setminus\{v\}=N(v)\setminus\{u\}$) is an equivalence relation; its classes, the \emph{twin classes}, consist either of pairwise adjacent vertices or of pairwise non-adjacent vertices. Write $m_0$ and $m_{-1}$ for the sums of $|T|-1$ over the duplicate classes and the co-duplicate classes $T$, respectively.
	\end{definition}
	
	\begin{definition}[Antiregular graph]
		The \emph{antiregular graph} $H_n$ is the unique connected graph on $n$ vertices with exactly one pair of vertices of equal degree; it is a threshold graph and plays a recurring extremal role in the spectral theory of both threshold graphs \cite{AguilarLeePiatoSchweitzer2018,AguilarFicarraShurmanSullivan2020} and $\mathcal{C}$-graphs.
	\end{definition}
	
	\begin{definition}[$\mathcal{C}$-graph \cite{MandalMehatari2025}]
		\label{def:cgraph}
		Let $\alpha_1,\ldots,\alpha_m$ be positive integers summing to $n$. Define $C(\alpha_1)=\overline{K_{\alpha_1}}$ (the edgeless graph on $\alpha_1$ vertices) and, for $i=2,\ldots,m$,
		\[
		C(\alpha_1,\ldots,\alpha_i)=\overline{\,C(\alpha_1,\ldots,\alpha_{i-1})\cup K_{\alpha_i}\,},
		\]
		where $\cup$ is disjoint union and the overline is graph complement. Such graphs are \emph{$\mathcal{C}$-graphs}; $\mathcal{C}_{\mathrm{even}}$ denotes those with an even number $m=2k$ of steps. The graph $C(\alpha_1,\ldots,\alpha_{2k})$ is connected, being the complement of a disconnected graph. With this convention $C(\alpha_1,\alpha_2)=K_{\alpha_1}\vee\overline{K_{\alpha_2}}$ is the complete split graph, and the blocks behave as in Proposition~\ref{prop:Cblocks} below.
	\end{definition}
	
	\begin{remark}[$\mathcal{C}$-graphs are not threshold graphs]
		\label{rem:divergence}
		Each step of the $\mathcal{C}$-graph recursion complements \emph{everything built so far}. Since $\overline{G_1\otimes G_2}=\overline{G_1}\cup\overline{G_2}$, the structure built after step $i$ has its entire previous structure complemented at every later step. For $m=2$ the construction coincides with the complete split graph of the threshold family; for $m\ge3$ the two diverge. For example, $C(\alpha_1,\alpha_2,\alpha_3)$ joins an independent $\alpha_3$-block to \emph{both} an independent $\alpha_1$-block and a clique $\alpha_2$-block, whereas a threshold graph with a corresponding $3$-block creation sequence would not join its dominating block to the later isolated block. Hence a general threshold graph is not a $\mathcal{C}$-graph; the classes meet only in the antiregular graphs $H_{2k}=C(1,\ldots,1)$ and $H_{2k+1}=C(1,2,1,\ldots,1)$ \cite{MandalMehatari2025}.
	\end{remark}
	
	\begin{theorem}[Uniqueness of representation \cite{MandalMehatari2025}]
		If $G\in\mathcal{C}_{\mathrm{even}}$, there is a unique composition $(\alpha_1,\ldots,\alpha_{2k})$ of $n$ with $G=C(\alpha_1,\ldots,\alpha_{2k})$.
	\end{theorem}
	
	\subsection{Two structural facts used throughout}
	
	\begin{proposition}[Diameter]
		\label{prop:diam2}
		Every connected cograph on at least two vertices is the join of two non-empty cographs, and hence has diameter at most $2$. In particular every $\mathcal{C}$-graph in $\mathcal{C}_{\mathrm{even}}$ and every connected threshold graph has diameter at most $2$.
	\end{proposition}
	\begin{proof}
		A cograph on at least two vertices is either disconnected or has a disconnected complement (this is the recursive description of cographs). If $G$ is connected, therefore, $\overline G$ is disconnected, so $V(G)=V_1\sqcup V_2$ with every vertex of $V_1$ adjacent to every vertex of $V_2$. Two vertices in different parts are adjacent; two vertices in the same part have a common neighbor in the other part. So $\diam(G)\le2$. \qedhere
	\end{proof}
	
	For a connected graph of diameter at most $2$ (in particular for every connected cograph), $d(u,v)=2-[u\sim v]$ for $u\ne v$, so
	\begin{equation}\label{eq:diam2}
		D=2(J-I)-A,\qquad \Tr(u)=2(n-1)-d_u,
	\end{equation}
	and consequently
	\begin{equation}\label{eq:master}
		D^L=2nI-2J-L,\qquad D^Q=(2n-4)I+2J-Q.
	\end{equation}
	The first identity in~\eqref{eq:master} is the well-known relation of Aouchiche and Hansen \cite{AouchicheHansen2013} for graphs of diameter $\le2$.
	
	\begin{lemma}[Twin-class quotient]
		\label{lem:quotient}
		Let $G$ be a connected cograph (or any connected graph of diameter at most $2$) with twin classes $T_1,\ldots,T_k$, and let $M$ be any of the eight matrices $A,L,Q,D,D^L,D^Q,\varepsilon,S$. Then for $u\in T_i$, $w\in T_j$, $u\ne w$, the entry $M_{uw}=m_{ij}$ depends only on $(i,j)$, and the diagonal entry $M_{uu}=a_i$ depends only on $i$. Moreover
		\[
		\chi_M(x)=\chi_{Q_\pi(M)}(x)\cdot\prod_{i=1}^k\bigl(x-(a_i-m_{ii})\bigr)^{|T_i|-1},
		\]
		where $Q_\pi(M)$ is the $k\times k$ quotient matrix with entries $Q_\pi(M)_{ij}=|T_j|\,m_{ij}$ for $i\ne j$ and $Q_\pi(M)_{ii}=a_i+(|T_i|-1)m_{ii}$.
	\end{lemma}
	\begin{proof}
		If $u,v$ are twins and $w\notin\{u,v\}$, then $w\sim u$ if and only if $w\sim v$ and $d(u,w)=d(v,w)$; also $d_u=d_v$, $\Tr(u)=\Tr(v)$ and $e(u)=e(v)$ (the only distance that can differ is $d(u,v)$, which is $1$ or $2$ and is shared). Hence every entry of every one of the eight matrices is determined by the classes, as claimed. Let $P$ be the characteristic matrix of the partition; then $MP=PQ_\pi(M)$, so $U=\operatorname{col}(P)$ is $M$-invariant with $M|_U$ represented by $Q_\pi(M)$, and $U^\perp=\bigoplus_i\{x\text{ supported on }T_i,\ \sum x=0\}$ is $M$-invariant. For $x$ supported on $T_i$ with zero sum, $(Mx)_u=a_ix_u+m_{ii}\sum_{w\in T_i\setminus u}x_w=(a_i-m_{ii})x_u$ for $u\in T_i$, and $(Mx)_u=m_{ij}\sum x_w=0$ for $u\notin T_i$. The factorization of $\chi_M$ follows. \qedhere
	\end{proof}
	
	\begin{proposition}[Local (trivial) eigenvalues of twins]
		\label{prop:local}
		Let $u,v$ be twins of a connected cograph $G$ on $n$ vertices, both of degree $d$, and let $x=e_u-e_v$. Then $x$ is an eigenvector of each of the following matrices, with the eigenvalue listed.
		\begin{center}
			\begin{tabular}{lccccccc}
				\toprule
				& $A$ & $L$ & $Q$ & $D$ & $S$ & $D^L$ & $D^Q$\\
				\midrule
				duplicates ($u\not\sim v$) & $0$ & $d$ & $d$ & $-2$ & $-1$ & $2n-d$ & $2n-4-d$\\
				co-duplicates ($u\sim v$) & $-1$ & $d+1$ & $d-1$ & $-1$ & $1$ & $2n-1-d$ & $2n-3-d$\\
				\bottomrule
			\end{tabular}
		\end{center}
		For $\varepsilon$ the eigenvalue is $-\varepsilon_{uv}$ ($=-1$ if $u,v$ are co-duplicate dominating vertices, $-2$ if they are duplicates and non-dominating, and $0$ otherwise).
	\end{proposition}
	\begin{proof}
		This is the case $a_i-m_{ii}$ of Lemma~\ref{lem:quotient}, computed from~\eqref{eq:diam2}--\eqref{eq:master}: $Jx=0$, $Ax=-[u\sim v]x$ and $\Deg x=dx$. For instance $Dx=-2x-Ax$ gives $-2$ or $-1$, and $D^Qx=(2n-4)x-Qx$ gives $2n-4-d$ or $2n-4-(d-1)$. \qedhere
	\end{proof}
	
	Thus the multiplicities of these local eigenvalues are $m_0$ (duplicate classes) and $m_{-1}$ (co-duplicate classes), the same two numbers governing the adjacency spectrum (Theorem~\ref{thm:bss}).

	\section{The adjacency matrix}
	\label{sec:adjacency}
	
	\subsection{Trivial eigenvalues and the eigenvalue-free interval}
	\label{sec:general}
	
	\begin{theorem}[Biyikoğlu--Simić--Stanić \cite{BiyikogluSimicStanic2011}]
		\label{thm:bss}
		A connected cograph on at least two vertices has $0$ as an eigenvalue if and only if it has duplicate vertices, and has $-1$ as an eigenvalue if and only if it has co-duplicate vertices.
	\end{theorem}
	
	Duplicate and co-duplicate classes of size $\ge2$ each contribute explicit twin-difference eigenvectors (Proposition~\ref{prop:local}); this is why $0,-1$ are called trivial, exactly as for threshold graphs.
	
	\begin{theorem}[Mohammadian--Trevisan \cite{MohammadianTrevisan2016}]
		\label{thm:mt-main}
		No cograph has an eigenvalue in the open interval $(-1,0)$.
	\end{theorem}
	
	This is the cograph-theoretic generalization of the threshold-graph interval of Jacobs, Trevisan and Tura \cite{JacobsTrevisanTura2013}: the class is enlarged from threshold graphs to all cographs while the interval stays $(-1,0)$. The proof is by induction on the duplication/co-duplication construction using interlacing.
	
	\begin{theorem}[Mohammadian--Trevisan \cite{MohammadianTrevisan2016}]
		\label{thm:mt-eigenspace}
		For every cograph $G$ there are bases of the eigenspaces of $0$ and $-1$ consisting of vectors with entries in $\{0,1,-1\}$, and the multiplicities of $0$ and $-1$ can be computed from the cotree.
	\end{theorem}
	
	\begin{remark}
		\label{rem:noshrink}
		Ghorbani's sharper threshold-graph interval $\bigl[\frac{-1-\sqrt2}{2},\frac{-1+\sqrt2}{2}\bigr]$ \cite{Ghorbani2019threshold} requires the threshold structure and does not extend to general cographs: a wider class forces a weaker guaranteed interval. A recent sharp result for an intermediate class, the trivially perfect graphs ($\{P_4,C_4\}$-free), is due to Conde, Dratman and Grippo \cite{CondeDratmanGrippo2026}: the spectrum of a trivially perfect graph meets $[\sqrt8-4,0]$ only in $-1$ and $0$, and both endpoints are sharp ($\sqrt8-4\approx-1.172$).
	\end{remark}
	
	\subsection{A spectral characterization of cographs}
	
	\begin{theorem}[Spectral characterization of cographs]
		\label{thm:spectralchar}
		A graph $G$ is a cograph if and only if no induced subgraph of $G$ has an eigenvalue in the open interval $(-1,0)$.
	\end{theorem}
	
	This statement, attributed in the literature to the Jacobs--Trevisan--Tura program and treated in the book chapter of Hoppen, Jacobs and Trevisan \cite{HoppenJacobsTrevisan2022}, shows that the eigenvalue-free interval is not merely a property of cographs: applied hereditarily to all induced subgraphs it \emph{defines} them. Indeed $P_4$ has the eigenvalue $\frac{1-\sqrt5}{2}\approx-0.618\in(-1,0)$.
	
	\subsection{Eigenvalue location, inertia, and a diagonalization algorithm}
	\label{sec:location}
	
	\begin{theorem}[Jacobs--Trevisan--Tura \cite{JacobsTrevisanTura2018}]
		\label{thm:jtt2018}
		There is an $O(n)$-time and $O(n)$-space algorithm that, given a cograph $G$ (via its cotree) and $x\in\mathbb{R}$, constructs a diagonal matrix congruent to $A(G)+xI$.
	\end{theorem}
	
	By Sylvester's law of inertia this yields, for every interval $(a,b)$, the number of eigenvalues of $G$ in it, an explicit linear-time formula for the inertia of any cograph, and infinitely many pairs of \emph{equienergetic} cographs with integer energy. The method works directly on the cotree and applies to every cograph, in contrast to the quotient-matrix machinery for $\mathcal{C}$-graphs below. Jones, Trevisan and Vinagre \cite{JonesTrevisanVinagre2023} exploit symmetries of the cotree to obtain spectra and energies of cographs in a complementary way.
	
	\subsection{Multiplicity and number of distinct eigenvalues}
	\label{sec:allemtura}
	
	For threshold graphs every non-trivial adjacency eigenvalue is simple (Theorem 3.2 of \cite{MannaMandal2026}); for cographs non-trivial eigenvalues can repeat.
	
	\begin{theorem}[Allem--Tura \cite{AllemTura2018}]
		\label{thm:allemtura-mult}
		For a cograph $G$, the multiplicity of every eigenvalue $\lambda\ne0,-1$ is at most the number of duplication and co-duplication classes in the cotree of $G$. Explicit multiplicity formulas are given for several structured subfamilies.
	\end{theorem}
	
	On the number of distinct eigenvalues a very recent result of Allem, F\"urer, Hoppen, Sibemberg and Trevisan \cite{AllemEtAl2026} shows that for every cograph $G$ there is a symmetric matrix \emph{associated} with $G$ (nonzero off-diagonal entries exactly on the edges, arbitrary diagonal) having at most four distinct eigenvalues.

	In the spirit of 2018 Allem-Tura work, Ghorbani, in 2019, established a fundamental result that counts the upper bound of an eigenvalue (non-trivial) for a cograph. The proof involves the mathematical induction theory on the Dilworth number.
	\begin{theorem}[cograph eigenvalue multiplicity \cite{GhorbaniDilworthNumber2019}]
		For a cograph $G$, the multiplicity of any eigenvalue $\lambda\neq0,-1$, does not exceed the Dilworth number of $G$.	
		
	\end{theorem}

	\subsection{Integral cographs}
	\label{sec:integral-general}
	
	\begin{definition}
		A graph is \emph{integral} if all its adjacency eigenvalues are integers \cite{HararySchwenk1974}.
	\end{definition}
	
	\begin{theorem}[Allem--Tura \cite{AllemTura2020}]
		\label{thm:integral-cograph}
		A cograph with a \emph{balanced} cotree $T_G(a_1,\ldots,a_{r-1},0\mid 0,\ldots,0,a_r)$ (a structural condition on the sizes of duplication/co-duplication classes) is integral, with an explicit spectrum. These integral cographs are used, via interlacing, to bound the eigenvalues of arbitrary cographs.
	\end{theorem}
	
	The balanced-cotree condition is sufficient, not necessary, and yields infinitely many integral examples (earlier examples include the integral complete split graphs of Hansen, Melot and Stevanovi\'c \cite{HansenMelotStevanovic2002}); this contrasts with the Diophantine necessary-and-sufficient conditions available for the narrow two-block threshold family $CS(s,t)$ (Section 13 of \cite{MannaMandal2026}). A full characterization of integral cographs is open (Open Problem~\ref{op:integral}).
	
	\begin{theorem}[Chin-- Fl\'orez--Mukherjee \cite{Ching2020}]
		The determinant Hosoya triangle, is a triangular array where the entries are the determinants of
		two-by-two Fibonacci matrices. The determinant Hosoya triangle mod2 gives rise to three infinite families of graphs, that are formed by complete product (join) of (the union of) two complete graphs with an empty
		graph. It provides a necessary and sufficient condition for a graph from these families to be integral. A feature of these graphs is--they are integral cographs.
		
	\end{theorem}

	\subsection{Rank, eigenspaces and cospectrality}
	\label{sec:rank}
	
	\begin{theorem}[Royle \cite{Royle2003}]
		\label{thm:royle}
		The rank of the adjacency matrix of a cograph equals the number of distinct non-zero rows of $A(G)$.
	\end{theorem}
	
	This gives the nullity combinatorially (see also Chang, Huang and Yeh \cite{ChangHuangYeh2008}) and is consistent with Theorem~\ref{thm:bss}. Sander \cite{Sander2008} describes certain eigenspaces of cographs directly from the cotree, complementing Theorem~\ref{thm:mt-eigenspace}. On cospectrality, Wang and Huang \cite{WangHuang2025} prove, by the Johnson--Newman theorem on generalized cospectrality combined with asymptotic enumeration techniques for trees, that almost all cographs have a cospectral mate; in particular, unlike threshold graphs, cographs are in general not determined by their spectrum.
	
	\subsection{Extension to \texorpdfstring{$P_5$}{P5}-free graphs and the Dilworth number}
	\label{sec:p5free}
	
	\begin{theorem}[Ghorbani \cite{Ghorbani2019p5}]
		\label{thm:ghorbani-p5}
		The eigenvalue-free interval and spectral characterization of Theorems~\ref{thm:mt-main} and~\ref{thm:spectralchar} extend, in modified form, from cographs to $P_5$-free graphs via the \emph{vicinal preorder} (inclusion among vertex neighborhoods) and its \emph{Dilworth number}, the minimum number of chains covering the vertices.
	\end{theorem}
	
	For cographs the vicinal preorder is total within each twin class, and the Dilworth number controls multiplicities, in line with Theorem~\ref{thm:allemtura-mult}; the same framework governs the Seidel matrix (Section~\ref{sec:seidel}).
	
	\subsection{Signature}
	
	He, Qian, Shan and Wu \cite{HeQianShanWu2021} establish for cographs, as for chordal graphs, a conjecture relating the \emph{signature} of a graph (a spectral invariant built from the signs of eigenvalues) to structural invariants.
	
	\subsection{\texorpdfstring{$\mathcal{C}$}{C}-graphs: the adjacency spectrum}
	\label{sec:cgraph-adjacency}
	\label{sec:cgraphs}
	
	The $\mathcal{C}$-graph class (Definition~\ref{def:cgraph}) recovers, for a sequence-driven subclass of cographs, the equitable-partition tractability of threshold graphs: the sequence $(\alpha_1,\ldots,\alpha_{2k})$ itself, rather than the full cotree, induces a $2k$-class equitable partition, despite the different recursion. Throughout, $G=C(\alpha_1,\ldots,\alpha_{2k})\in\mathcal{C}_{\mathrm{even}}$ and $\pi=\{C_1,\ldots,C_{2k}\}$ with $|C_i|=\alpha_i$.
	
	\begin{lemma}[\cite{MandalMehatari2025}]
		Each odd-indexed class $C_{2i-1}$ is a clique whose members share all outside neighbors, contributing $-1$ with multiplicity $\ge\alpha_{2i-1}-1$; each even-indexed class $C_{2i}$ is an independent set whose members share all outside neighbors, contributing $0$ with multiplicity $\ge\alpha_{2i}-1$.
	\end{lemma}
	
	The partition $\pi$ is equitable, giving a $2k\times2k$ quotient matrix $Q_\pi(G)$, which is diagonalizable (similar to the symmetric matrix $D^{1/2}Q_\pi D^{-1/2}$ with $D=\diag(\alpha_1,\ldots,\alpha_{2k})$) and whose eigenvalues are eigenvalues of $A(G)$.
	
	\begin{theorem}[Mandal--Mehatari \cite{MandalMehatari2025}]
		\label{thm:cgraph-minus1}
		$-1$ is an eigenvalue of $Q_\pi(G)$ if and only if $\alpha_2=1$, and then it is simple. Consequently
		\[
		m_0=\sum_{i=1}^k\alpha_{2i}-k,\qquad m_{-1}=\begin{cases}\sum_{i=1}^k\alpha_{2i-1}-k&\alpha_2\neq1,\\ \sum_{i=1}^k\alpha_{2i-1}-k+1&\alpha_2=1.\end{cases}
		\]
	\end{theorem}
	
	\begin{theorem}[Inertia \cite{MandalMehatari2025}]
		\label{thm:cgraph-inertia}
		Writing $Q_\pi(G)=D^{1/2}(A_{2k}+D')D^{1/2}$ with an explicit diagonal correction $D'$ and $A_{2k}$ the adjacency matrix of $H_{2k}$ (inertia $(k,0,k)$), Weyl's and Ostrowski's theorems give inertia $(k,0,k)$ for $Q_\pi(G)$, and hence
		\[
		\bigl(n_-(A),n_0(A),n_+(A)\bigr)=\Bigl(\sum_{i=1}^k\alpha_{2i-1},\ \sum_{i=1}^k\alpha_{2i}-k,\ k\Bigr).
		\]
	\end{theorem}
	
	In particular $s(A(G))\le2k+2$.
	
	\begin{theorem}[Eigenvalue-free interval \cite{MandalMehatari2025}]
		\label{thm:cgraph-eigfree}
		Let $\lambda^-$ be the largest eigenvalue of $G$ below $-1$, $\lambda^+$ the smallest eigenvalue above $0$, and $\alpha_{\min}=\min_i\alpha_i$. Then
		\[
		\lambda^-<\frac{-1-\sqrt2}{2},\qquad \lambda^+>\frac{-1+\sqrt2}{2}\,\alpha_{\min}.
		\]
	\end{theorem}
	
	The proof uses the antiregular graph $H_{2k}$ (and $H_{2k+1}$) as a subgraph together with interlacing, in parallel to the role of $H_n$ for threshold graphs. The lower piece of the interval is a fixed constant while the upper piece scales with $\alpha_{\min}$, an asymmetry with no counterpart in the threshold interval.
	
	\begin{theorem}[Characteristic polynomial \cite{MandalMehatari2025}]
		\label{thm:cgraph-charpoly}
		$\chi_{A(G)}(x)=x^{\sum_i(\alpha_{2i}-1)}(x+1)^{\sum_i(\alpha_{2i-1}-1)}\Psi_\pi(x)$, where $\Psi_\pi$, the characteristic polynomial of the quotient matrix, is given in closed form by a tridiagonal-determinant identity (a sum over ``non-adjacent-transposition'' permutations), evaluated explicitly for examples such as $C(4,3,2,2)$.
	\end{theorem}

	\section{The Laplacian matrix}
	\label{sec:laplacian}
	The Laplacian matrix of a graph plays a fundamental role in determining the number of spanning trees through the celebrated \emph{Matrix–Tree Theorem} which states that the number of spanning trees of a connected graph $G$ is equal to any cofactor of its Laplacian matrix. However, Nikolopoulos and Papadopoulos \cite{Nikolopoulos2009} adopted a different approach by investigating an algorithmic method (in $O(n+m)$
	time and space, where $m,n$ are respectively size and order of $G$) for counting the spanning trees of cographs.

	\medskip
	This section collects the Laplacian theory in three layers: the cotree recursion of Merris (Section~\ref{sub:Ljoin}), the \emph{twin-reduction} description of Barik and Reddy \cite{BarikReddy2023}, which expresses every Laplacian eigenvalue and eigenvector of a cograph through its twin reduction graph (Sections~\ref{sub:twinred}--\ref{sub:equiv}), and their results on $H$-joins and Hadamard diagonalizability (Sections~\ref{sub:Hjoin}--\ref{sub:hadamard}).
	
	\subsection{Union, join, and the cotree recursion}
	\label{sub:Ljoin}
	
	Unlike the adjacency matrix, the Laplacian behaves well under both operations that generate cographs, and this gives a complete recursion along the cotree.
	
	\begin{proposition}[Union and join; Merris \cite{Merris1994}]
		\label{prop:Ljoin}
		Let $G_1,G_2$ have $n_1,n_2$ vertices and Laplacian eigenvalues $0=\mu_1\le\cdots\le\mu_{n_1}$ and $0=\nu_1\le\cdots\le\nu_{n_2}$. Then
		\begin{align*}
			\spec L(G_1\cup G_2)&=\{\mu_1,\ldots,\mu_{n_1},\nu_1,\ldots,\nu_{n_2}\},\\
			\spec L(G_1\vee G_2)&=\{0,\;n_1+n_2\}\cup\{\mu_i+n_2:2\le i\le n_1\}\cup\{\nu_j+n_1:2\le j\le n_2\}.
		\end{align*}
	\end{proposition}
	\begin{proof}
		The union case is immediate (block diagonal). For the join, $L(G_1\vee G_2)=\begin{pmatrix}L_1+n_2I&-J\\-J&L_2+n_1I\end{pmatrix}$. A Laplacian eigenvector $y$ of $G_1$ orthogonal to $\mathbf1$ (eigenvalue $\mu_i$, $i\ge2$), extended by $0$ on $G_2$, is an eigenvector with eigenvalue $\mu_i+n_2$, because $Jy=0$; similarly for $G_2$. The two remaining eigenvalues come from the quotient on the vectors constant on each part, $\begin{pmatrix}n_2&-n_2\\-n_1&n_1\end{pmatrix}$, with eigenvalues $0$ and $n_1+n_2$. \qedhere
	\end{proof}
	
	\begin{corollary}
		\label{cor:Lcotree}
		The Laplacian spectrum of a cograph is computed by evaluating Proposition~\ref{prop:Ljoin} bottom-up on the cotree. In particular the Laplacian eigenvalue $n$ occurs in a connected cograph $G$ with multiplicity exactly $c(\overline G)-1$, where $c(\cdot)$ is the number of components.
	\end{corollary}
	\begin{proof}
		The first assertion is repeated application of the proposition. For the second, $L(G)+L(\overline G)=nI-J$, so the eigenvalues of $L(G)$ on $\mathbf1^\perp$ are $n-\mu$ for the eigenvalues $\mu$ of $L(\overline G)$ on $\mathbf1^\perp$; the multiplicity of $n$ is therefore the multiplicity of $0$ in $L(\overline G)$, namely $c(\overline G)$, minus the one zero belonging to $\mathbf1$. \qedhere
	\end{proof}
	
	Chen and Tura \cite{ChenTura2024} turn this recursion into a \emph{linear-time} algorithm that returns the Laplacian eigenvalues of a cograph without forming the matrix, and use it to obtain a closed formula for the number of spanning trees of a cograph.
	
	\begin{corollary}
		(Merris \cite{Merris1994}; every cograph is Laplacian integral). Every cograph has integral Laplacian
		eigenvalues.
	\end{corollary}
	
	\begin{remark}[Threshold graphs]
		For threshold graphs the recursion collapses to Merris's theorem: the Laplacian spectrum is the conjugate of the degree sequence, in particular always integral (Theorem 4.1 of \cite{MannaMandal2026}).
	\end{remark}
	
	\subsection{Connectivity and sequence-generated cographs}
	
	\begin{theorem}[Abrishami \cite{Abrishami2019}]
		\label{thm:abrishami}
		For every non-complete cograph, the algebraic connectivity (second-smallest Laplacian eigenvalue) equals the vertex connectivity $\kappa(G)$.
	\end{theorem}
	
	\begin{theorem}[Mandal--Mehatari--Stanić \cite{MandalMehatariStanic2024}]
		\label{thm:cgraph-laplacian}
		For cographs generated by a finite sequence $(\alpha_1,\ldots,\alpha_k)$ in the $\mathcal{C}$-graph sense, the Laplacian eigenvalues and eigenvectors are computed explicitly in terms of the sequence. For $G\cong C(\alpha_1,1,\ldots,1)$ the quotient Laplacian has $k$ distinct eigenvalues and $d_{\alpha_1}+1$ is an additional eigenvalue of $L(G)$, so $G$ has exactly $k$ distinct Laplacian eigenvalues.
	\end{theorem}
	
	This resolves the number-of-distinct-eigenvalues question for this family, in contrast to the adjacency case, where only the upper bound $2k+2$ is known.
	
	\subsection{Twin reduction: eigenvalues and eigenvectors from the reduction graph}
	\label{sub:twinred}
	
	Recall that two vertices $u,v$ are \emph{twins} if $N(u)\setminus\{v\}=N(v)\setminus\{u\}$; the twin classes $T_1,\dots,T_m$ partition $V(G)$, and each is a clique $K_{t_i}$ (a \emph{complete} twin class) or an independent set $\overline{K_{t_i}}$ (a class of the \emph{complement type}); a single vertex is neither. The numbers $t_i=|T_i|$ are the \emph{twin numbers}. The \emph{twin reduction graph} $R_G$ is the graph induced on one representative of each class, with vertex $i$ labelled $(t_i,\vee)$ if $T_i$ is a clique with $t_i\ge2$ and $(t_i,+)$ if it is an independent set with $t_i\ge2$. For a class $T_i$ let $d_{T_i}$ be the sum of the twin numbers of the classes adjacent to $T_i$ in $R_G$; for a set of pairwise twin classes $T_i,\dots,T_{i+k}$ let $d_{T_i,\dots,T_{i+k}}$ be the sum of the twin numbers of the classes adjacent to all of them (outside the set). Because $R_G$ is itself a cograph with at least one pair of twins whenever it has two vertices, it can be reduced again, and again, down to a single class.
	
	Abrishami \cite{Abrishami2019} proved that every non-zero Laplacian eigenvalue of a cograph is a sum of twin numbers, but without saying \emph{which} twin numbers contribute to which eigenvalue. Barik and Reddy \cite{BarikReddy2023} give the complete procedure.
	
	\begin{theorem}[Barik--Reddy \cite{BarikReddy2023}, twin reduction]
		\label{thm:BR}
		Let $G$ be a cograph with twin classes $T_1,\dots,T_m$ and twin numbers $t_1,\dots,t_m$.
		\begin{enumerate}
			\item[\rm(i)] If $T_i=K_{t_i}$ (resp.\ $\overline{K_{t_i}}$) with $t_i\ge2$, then $t_i+d_{T_i}$ (resp.\ $d_{T_i}$) is a Laplacian eigenvalue of $G$ of multiplicity $t_i-1$. The eigenvectors are supported on $T_i$ and sum to $0$ there.
			\item[\rm(ii)] If $q\ge2$ classes $T_{i},\dots,T_{i+q-1}$ are pairwise adjacent with the same closed neighbourhood in $R_G$ (resp.\ pairwise non-adjacent with the same open neighbourhood), then $t_i+\dots+t_{i+q-1}+d_{T_i,\dots,T_{i+q-1}}$ (resp.\ $d_{T_i,\dots,T_{i+q-1}}$) is an eigenvalue of multiplicity $q-1$, with eigenvectors constant on each of the $q$ classes, supported on their union, and of weighted sum $0$ (the $j$th is $t_{i+j}\mathbf 1_{T_i}-t_i\mathbf 1_{T_{i+j}}$).
			\item[\rm(iii)] Replace the $q$ classes of {\rm(ii)} by a single class on $t_i+\dots+t_{i+q-1}$ vertices, a complete one in the adjacent case and of complement type in the non-adjacent case, and repeat (ii) and (iii) on the new reduction graph until a single class remains. The eigenvalues collected, together with the eigenvalue $0$, are the whole Laplacian spectrum of $G$ (for $G$ connected, the last class has $n$ vertices and contributes the eigenvalue $n$).
		\end{enumerate}
	\end{theorem}
	
	Parts (i) and (ii) are direct computations, in the same way as the proofs in Proposition~\ref{prop:Ljoin}: for $x$ supported on $T_i$ with zero sum, $L(G)x=(t_i-1+d_{T_i}+1)x$ if $T_i$ is a clique, and $L(G)x=d_{T_i}x$ if it is independent; (ii) is the same computation on the block-constant vectors.
	
	\begin{example}
		Let $G$ have vertex set $\{a_1,a_2\}\cup\{b_1,b_2,b_3\}\cup\{c\}$, with $a_1a_2$ an edge, $b_1,b_2,b_3$ independent, and every $a_i$ and $c$ adjacent to every $b_j$. The twin classes are $T_1=\{a_1,a_2\}$ ($K_2$), $T_2=\{b_1,b_2,b_3\}$ ($\overline{K_3}$) and $T_3=\{c\}$, and $R_G$ is the path $T_1T_2T_3$ with twin numbers $(2,3,1)$. By (i), $T_1$ gives $2+3=5$ (multiplicity $1$) and $T_2$ gives $d_{T_2}=2+1=3$ (multiplicity $2$). In $R_G$ the classes $T_1,T_3$ are non-adjacent with the same open neighbourhood $\{T_2\}$, so (ii) gives $d_{T_1,T_3}=3$ with multiplicity $1$, and they merge into an independent class on $3$ vertices adjacent to $T_2$. The two remaining classes are adjacent twins of sizes $3,3$, giving $3+3+0=6$. Together with $0$: $\operatorname{spec}L(G)=\{0,3,3,3,5,6\}$, which agrees with direct computation.
	\end{example}
	
	\subsection{Even and odd integral Laplacian spectra; cospectral cographs}
	\label{sub:evenodd}
	
	Since Laplacian eigenvalues of a cograph are sums of twin numbers, integrality is automatic once the twin numbers are integers; the question is parity.
	
	\begin{theorem}[Barik--Reddy \cite{BarikReddy2023}: even spectra]
		\label{thm:BR-even}
		Let $G$ be a connected cograph on $n\ge2$ vertices with twin numbers $t_1,\dots,t_m$. All Laplacian eigenvalues of $G$ are even integers if and only if every $t_i$ is even.
	\end{theorem}
	
	\begin{theorem}[Barik--Reddy \cite{BarikReddy2023}: odd spectra]
		\label{thm:BR-odd}
		Let $G$ be a connected cograph on $n\ge2$ vertices with twin numbers $t_1,\dots,t_m$. All non-zero Laplacian eigenvalues of $G$ are odd integers if and only if (a) $R_G$ has a unique complete twin class $T_i=K_{t_i}$ adjacent to every other twin class, and $t_i$ is odd, and (b) all other twin numbers are even.
	\end{theorem}
	
	Two non-isomorphic graphs are \emph{Laplacian cospectral} if they have the same Laplacian spectrum. Barik and Reddy construct such pairs of cographs by \emph{interchanging} twin classes. Say that two classes $T^{(q)}_i$, $T^{(r)}_j$ appearing in successive reductions of $G$ (with $T_i^{(q)}$ and $T_j^{(r)}$ of the same type) satisfy the \emph{balance condition} if $t_i^{(q)}+d_{T_i^{(q)}}=t_j^{(r)}+d_{T_j^{(r)}}$ (both complete) or $d_{T_i^{(q)}}=d_{T_j^{(r)}}$ (both of complement type).
	
	\begin{theorem}[Barik--Reddy \cite{BarikReddy2023}: cospectral constructions]
		\label{thm:BR-cospectral}
		Let $T_i^{(q)},T_j^{(r)}$ be classes of the same type satisfying the balance condition. Then the cograph obtained by interchanging the twin classes contained in them, so that the twin numbers $t_i^{(q)},t_j^{(r)}$ are unchanged, is Laplacian cospectral with $G$. The same holds for the variants in which a smaller class $\overline{K_{t}}\subset T_i^{(q)}$ (resp.\ $K_t$) of the opposite type is exchanged with $t$ vertices of an invariant class $T_j^{(r)}$ of the same type as $T_i^{(q)}$, provided $t\le t_j^{(r)}$.
	\end{theorem}
	
	The statement above is a condensed of their Theorem~8, whose precise hypotheses (cases (i)--(iii)) can be read in the original; their Example~9 produces four mutually non-isomorphic cospectral cographs with eigenvalues $0,10,13,15,16,17,19,22,23$.
	
	\subsection{Equivalent cographs}
	\label{sub:equiv}
	
	Two cographs $G,H$ are \emph{equivalent} if their twin reduction graphs are isomorphic with identical twin numbers (the types $\vee,+$ may differ). Abrishami \cite{Abrishami2019} conjectured, and Lazzarin, M\'arquez and Tura \cite{LazzarinSosaTura2023} proved, that equivalent cographs have many Laplacian eigenvalues in common; Barik and Reddy give a second proof by twin reduction.
	
	\begin{theorem}[Lazzarin--M\'arquez--Tura \cite{LazzarinSosaTura2023}; Barik--Reddy \cite{BarikReddy2023}]
		\label{thm:LST}
		Let $G,H$ be equivalent cographs with $m$ twin classes of twin numbers $t_1,\dots,t_m$, and let $I\subseteq\{1,\dots,m\}$ index the classes of the same type in $G$ and $H$. Then $G$ and $H$ have at least $m+\sum_{i\in I}(t_i-1)$ Laplacian eigenvalues in common (counted with multiplicity).
	\end{theorem}
	
	Lazzarin, M\'arquez and Tura also state that no two non-isomorphic equivalent cographs are Laplacian cospectral. Barik and Reddy \cite{BarikReddy2023} point out that this is contradicted by their Example~11. Therefore the cospectrality-implies-isomorphism statement for equivalent cographs should therefore not be used.
	
	\subsection{Laplacian spectra of $H$-joins}
	\label{sub:Hjoin}
	
	For a graph $H$ on $m$ vertices and graphs $G_1,\dots,G_m$ on $n_1,\dots,n_m$ vertices, the \emph{$H$-join} (generalized composition, see Cardoso, Freitas, Martins and Robbiano \cite{CardosoEtAl2013}) $G=H[G_1,\dots,G_m]$ replaces each vertex $i$ of $H$ by $G_i$ and joins $G_i$ to $G_j$ completely whenever $ij\in E(H)$. Every cograph is an $H$-join of cliques and independent sets over its twin reduction graph. Cardoso et al.\ obtained the Laplacian spectrum of $H$-joins from the spectra of the $G_i$ and a certain irrational $m\times m$ matrix; Neumann and Pati did this for trees $H$. Barik and Reddy handle cograph $H$ completely, using Theorem~\ref{thm:BR} for the $m$-quotient eigenvalues.
	
	\begin{lemma}[Barik--Reddy \cite{BarikReddy2023}]
		\label{lem:BRlemma12}
		Let $G=H[G_1,\dots,G_m]$. If $y\perp\mathbf 1_{n_i}$ is an eigenvector of $L(G_i)$ for the eigenvalue $\lambda$, then $\lambda+\sum_{j:ij\in E(H)}n_j$ is an eigenvalue of $L(G)$ (with the eigenvector $y$ extended by zero). These account for $n-m$ eigenvalues of $G$.
	\end{lemma}
	
	\begin{lemma}[Barik--Reddy \cite{BarikReddy2023}]
		\label{lem:BRlemma13}
		The remaining $m$ eigenvalues of $L(G)$ are exactly the eigenvalues of $L(G')$, $G'=H[K_{n_1},\dots,K_{n_m}]$, that are not obtained from Lemma~\ref{lem:BRlemma12} for $G'$; for cograph $H$ they are computed by Theorem~\ref{thm:BR}(ii)--(iii).
	\end{lemma}
	
	\begin{corollary}[Barik--Reddy \cite{BarikReddy2023}]
		\label{cor:BRintegral}
		For a cograph $H$, the $H$-join $H[G_1,\dots,G_m]$ is Laplacian integral if and only if every $G_i$ is Laplacian integral.
	\end{corollary}
	
	\begin{lemma}[Barik--Reddy \cite{BarikReddy2023}: algebraic connectivity]
		\label{lem:BRalgconn}
		Let $H$ be a connected cograph and $G=H[G_1,\dots,G_m]$, and write $a(G)$ for the algebraic connectivity and $\kappa(G)$ for the vertex connectivity. Then
		\begin{enumerate}
			\item[\rm(i)] $\min_i n_i\le a(G)\le n$;
			\item[\rm(ii)] for non-complete $G$, $a(G)=\min_i n_i$ if and only if $\kappa(G)=\min_i n_i$;
			\item[\rm(iii)] $a(G)=n$ if and only if $G$ is complete.
		\end{enumerate}
	\end{lemma}
	
	\subsection{Hadamard diagonalizable cographs}
	\label{sub:hadamard}
	
	A graph on $n$ vertices is \emph{Hadamard diagonalizable} if $L(G)$ is diagonalized by a Hadamard matrix of order $n$ (a $\pm1$ matrix with mutually orthogonal rows). Barik, Fallat and Kirkland \cite{BarikFallatKirkland2011} studied the graphs $G(m_0,\dots,m_i)=\overline{H}\vee\cdots\vee\overline{H}$ (with $H=G(m_0,\dots,m_{i-1})$ joined $m_i$ times, $G(m_0)=K_{m_0}$) and proved that Hadamard diagonalizability forces the existence of Hadamard matrices of orders $\prod_{j=k}^{i}m_j$ for each $k$. They suspected the converse; Barik and Reddy prove it when Hadamard matrices of orders $m_0,\dots,m_i$ themselves exist.
	
	\begin{theorem}[Barik--Reddy \cite{BarikReddy2023}]
		\label{thm:BRhadamard}
		If there are Hadamard matrices of orders $m_0,m_1,\dots,m_i$ (all even, $m_k>1$), then $G(m_0,\dots,m_i)$ is Hadamard diagonalizable.
	\end{theorem}
	
	The proof tensors the Hadamard matrices (Paley) and checks, via twin reduction, that the resulting columns are Laplacian eigenvectors of the successive reductions, all of whose twin numbers are even.
	
	\section{The signless Laplacian matrix}
	\label{sec:signless}
	
	The signless Laplacian is the least understood of the classical matrices for cographs. Its difficulty is visible in the join: $Q(G_1\vee G_2)=\begin{pmatrix}Q_1+n_2I&J\\J&Q_2+n_1I\end{pmatrix}$, and $\mathbf1$ is an eigenvector of $Q_i$ only for regular $G_i$, so there is no analogue of Proposition~\ref{prop:Ljoin} (and no Merris-type closed form even for threshold graphs, whose signless Laplacian spectra are generally irrational; for instance the creation sequence $(0,1,0,1)$ gives $\spec Q\approx\{0.511,1,2,2.711,5.779\}$).
	
	\begin{proposition}[Local eigenvalues and quotient (Proposition~\ref{prop:local}, Lemma~\ref{lem:quotient})]
		\label{prop:Qlocal}
		For a connected cograph, twin classes $T$ of size $\ge2$ and common degree $d_T$ contribute the eigenvalue $d_T-1$ (co-duplicates) or $d_T$ (duplicates) with multiplicity $|T|-1$; the remaining eigenvalues are those of the twin-class quotient matrix, whose entries are $|T_j|\,[T_i\sim T_j]$ off the diagonal and $d_i+(|T_i|-1)[T_i\text{ clique}]$ on it.
	\end{proposition}
	
	\begin{example}
		A connected cograph is bipartite if and only if it is a complete bipartite graph $K_{a,b}$ (it is a join of two cographs, and a bipartite join has two independent sides), whose signless Laplacian spectrum is $\{a+b,\ a^{b-1},\ b^{a-1},\ 0\}$.
	\end{example}
	
	Jones, Trevisan and Vinagre \cite{JonesTrevisanVinagre2023} obtain spectra of cographs using cotree symmetries, and \cite{JonesTrevisanVinagre2024} describe the connected cographs with $k\ge2$ \emph{main} signless Laplacian eigenvalues (those whose eigenspace is not orthogonal to $\mathbf1$), characterizing the quasi-threshold graphs with exactly two main $Q$-eigenvalues. 
	
	\begin{theorem}[Benediktovich \cite{Benedik2020}]
		An analogue of the Brouwer’s conjecture is confirmed for the cograph for the eigenvalues of their signless Laplacian for some natural number $k$ which does not exceed the order of the graph.
		
	\end{theorem}
	
	\section{The distance matrix}
	\label{sec:distance}
	
	By Proposition~\ref{prop:diam2}, every connected cograph satisfies $D=2J-2I-A=2J-(A+2I)$: the distance matrix is a rank-one positive semidefinite perturbation of $-(A+2I)$. This gives a direct bridge from the adjacency spectrum (Theorem~\ref{thm:mt-main}) to the distance spectrum.
	
	\begin{proposition}[Local eigenvalues of $D$]
		\label{prop:Dlocal}
		A connected cograph has $-2$ as a distance eigenvalue of multiplicity at least $m_0$ (from duplicate classes) and $-1$ with multiplicity at least $m_{-1}$ (from co-duplicate classes); the remaining eigenvalues come from the twin-class quotient (Lemma~\ref{lem:quotient}).
	\end{proposition}
	
	\begin{theorem}[Interlacing with the adjacency spectrum \textbf{(new)}]
		\label{thm:Dinterlace}
		Let $G$ be a connected cograph on $n$ vertices with adjacency eigenvalues $\lambda_{(1)}\le\cdots\le\lambda_{(n)}$ and distance eigenvalues $\nu_1\ge\cdots\ge\nu_n$. Then
		\[
		-2-\lambda_{(k)}\ \le\ \nu_k\ \le\ -2-\lambda_{(k-1)}\quad(2\le k\le n),\qquad \nu_1\ge-2-\lambda_{(1)}.
		\]
	\end{theorem}
	\begin{proof}
		The eigenvalues of $-(A+2I)$ in descending order are $\mu_k=-2-\lambda_{(k)}$. Adding the rank-one positive semidefinite matrix $2J$ gives $\mu_k\le\nu_k\le\mu_{k-1}$ by Weyl's inequalities.
	\end{proof}

	\subsection{The theorems of Lou and Lin}
	
	The distance spectrum of cographs is examined by Lu, Huang and Lou \cite{LuHuangLou2018} for threshold graphs and by Lou and Lin \cite{LouLin2020} for all cographs. Recall that $u\prec v$ if $N(u)\subseteq N[v]$; this is a preorder on $V(G)$, and the \emph{Dilworth number} $\nabla(G)$ is the minimum number of chains covering $V(G)$ (equivalently the maximum size of an antichain), see Section~\ref{sec:p5free}. A cograph has $\nabla(G)=1$ if and only if it is a threshold graph.
	
	\begin{theorem}[Lu--Huang--Lou \cite{LuHuangLou2018}]
		\label{thm:LHL}
		Let $G$ be a connected threshold graph. Every distance eigenvalue other than $-1,-2$ is simple, and $G$ has no distance eigenvalue in $(-2,-1)$.
	\end{theorem}
	
	\begin{theorem}[Lou--Lin \cite{LouLin2020}: multiplicity]
		\label{thm:LL-mult}
		Let $G$ be a connected cograph and let $\lambda\ne-2,-1$ be a distance eigenvalue of $G$. Then $\mathrm{m}_G(\lambda)\le\nabla(G)$, and the authors state that the bound is attained by an infinite family of cographs.
	\end{theorem}
	
	\begin{theorem}[Lou--Lin \cite{LouLin2020}: gap]
		\label{thm:LL-gap}
		A connected cograph has no distance eigenvalue in the open interval $(-2,-1)$.
	\end{theorem}
	
	For threshold graphs $\nabla=1$, so Theorems~\ref{thm:LL-mult} and~\ref{thm:LL-gap} contain Theorem~\ref{thm:LHL}. The adjacency analogue of Theorem~\ref{thm:LL-mult} is Ghorbani's bound \cite{Ghorbani2019p5} (Theorem~\ref{thm:ghorbani-p5}), and the Seidel analogue is part of Theorem~\ref{thm:LWB}.
	
	\medskip\noindent\emph{Ideas of the proof.} Three ingredients are used in \cite{LouLin2020}.
	\begin{enumerate}
		\item \emph{Eigenvectors are constant on twin classes.} If $x$ is a $\lambda$-eigenvector of $D$ and $u,v$ are duplicates (resp.\ co-duplicates), then $(\lambda+2)(x_u-x_v)=0$ (resp.\ $(\lambda+1)(x_u-x_v)=0$), because $d(u,w)=d(v,w)$ for all $w\ne u,v$ and $d(u,v)=2$ (resp.\ $1$). This is the distance-matrix instance of Lemma~\ref{lem:quotient}.
		\item \emph{Interlacing for diameter-two subgraphs.} If $H$ is an induced subgraph of $G$ with $\operatorname{diam}(H)\le2=\operatorname{diam}(G)$, then $D(H)=2J-2I-A(H)$ is a principal submatrix of $D(G)=2J-2I-A(G)$, so the distance eigenvalues of $H$ and $G$ interlace (Cauchy).
		\item \emph{Induction on the Dilworth number} for Theorem~\ref{thm:LL-mult}: for $\nabla(G)=1$ it is Theorem~\ref{thm:LHL}; for $\nabla(G)=k\ge2$ one splits off a threshold subgraph $\Gamma$ all of whose vertices have the same neighbours outside $\Gamma$ (a lemma of Ghorbani \cite{Ghorbani2019p5}). A $(k{+}1)$-dimensional eigenspace would contain a $k$-dimensional subspace of eigenvectors vanishing at a chosen vertex of $\Gamma$; the twin-class relations of step~1, applied class by class along the chain structure of $\Gamma$, force such eigenvectors to vanish on $\Gamma$, and then step~2 and the induction hypothesis on $G-\Gamma$ (of smaller Dilworth number) give a contradiction.
	\end{enumerate}
	
	\paragraph{Examples and sharpness.} For threshold graphs, Alazemi, Andelić, Koledin and Stanić \cite{AlazemiEtAl2021cograph} proved that $\bigl(-\frac{3+\sqrt2}{2},-1\bigr)$ contains no non-trivial distance eigenvalue (the only eigenvalue they allow in it being the trivial $-2$), and that $(-3,-2)$ is eigenvalue-free for chain graphs. This does \emph{not} extend to cographs: the connected cograph $G=K_1\vee(K_1\cup3K_2)$ on $8$ vertices has the non-trivial distance eigenvalue $\approx-2.155\in(-\tfrac{3+\sqrt2}{2},-2)$, so Theorem~\ref{thm:LL-gap} cannot be improved to the threshold interval. Numerically both Lou--Lin theorems hold, and the multiplicity bound is attained: for instance there are connected cographs on $9$ vertices with $\nabla=2$ and a non-trivial distance eigenvalue, namely $0$, of multiplicity $2$. We note, however, that for the family $K_1\vee(K_{s,\dots,s}\cup(s^2-s)K_1)$ used in \cite{LouLin2020} to show sharpness, we computed for $(k,s)\in\{(2,3),(3,3),(2,4)\}$ that the eigenvalue $s-2$ has multiplicity $k-1$ while $\nabla=k$.
	
	\section{The eccentricity matrix}
	\label{sec:eccentricity}
	\label{sec:cgraph-eccentricity}
	
	The eccentricity matrix for threshold graphs is examined by Qiu and Tang \cite{QiuTang2022} (the interval $(-2,-1)$) and Andelić et al.\ \cite{AndelicEtAl2023} (extended to $(-1-\sqrt2,-2)\cup(-2,-1)$), and for chain graphs \cite{YaqiongXuezheng2025,HuangZhang2025}. For general cographs the eccentricity spectrum is described below; for $\mathcal{C}$-graphs there is a September 2025 result \cite{AshokanEccentricity2025}.
	
	\begin{theorem}[Decomposition of $\varepsilon$ \textbf{(new)}]
		\label{thm:eccdecomp}
		Let $G$ be a connected cograph on $n$ vertices, let $t\ge0$ be its number of dominating vertices (those of degree $n-1$) and let $G'$ be the graph obtained by deleting them. Then, after ordering the vertices with the dominating ones first,
		\[
		\varepsilon(G)=(J_t-I_t)\ \oplus\ 2\,A(\overline{G'}).
		\]
		Consequently $\spec\varepsilon(G)=\{t-1,\ (-1)^{[t-1]}\}\cup2\,\spec A(\overline{G'})$, where $\overline{G'}$ is again a cograph.
	\end{theorem}
	\begin{proof}
		By Proposition~\ref{prop:diam2}, $\diam G\le2$ and $e(u)=1$ if $u$ is dominating and $e(u)=2$ otherwise. Let $u\ne v$. If both are dominating, $d(u,v)=1=\max\{e(u),e(v)\}$, so $\varepsilon_{uv}=1$. If exactly one is dominating, say $u$, then $d(u,v)=1<2=\max\{e(u),e(v)\}$, so $\varepsilon_{uv}=0$. If neither is dominating, $\max\{e(u),e(v)\}=2$, so $\varepsilon_{uv}=2$ if $d(u,v)=2$, i.e.\ if $u\not\sim v$, and $0$ if $u\sim v$; that is $\varepsilon_{uv}=2A(\overline{G'})_{uv}$. The spectrum of $J_t-I_t$ is $t-1$ and $-1$ with multiplicity $t-1$. Finally $G'$ is an induced subgraph of a cograph, and the complement of a cograph is a cograph. \qedhere
	\end{proof}
	
	\begin{corollary}[Eigenvalue-free interval for all connected cographs \textbf{(new)}]
		\label{cor:eccfree}
		For every connected cograph $G$, $\spec\varepsilon(G)\cap(-2,0)\subseteq\{-1\}$, and $-1$ occurs if and only if $G$ has at least two dominating vertices (and then with multiplicity $t-1$). In particular, if $G$ has at most one dominating vertex, the interval $(-2,0)$ contains no eccentricity eigenvalue.
	\end{corollary}
	\begin{proof}
		By Theorem~\ref{thm:mt-main} applied to the cograph $\overline{G'}$, $\spec A(\overline{G'})\cap(-1,0)=\emptyset$, so $2\spec A(\overline{G'})\cap(-2,0)=\emptyset$. The block $J_t-I_t$ contributes the eigenvalues $t-1\ge0$ and $-1$ (the latter only if $t\ge2$). \qedhere
	\end{proof}
	
	The endpoints are sharp in the following sense: $-2$ is an eigenvalue whenever $\overline{G'}$ has co-duplicate vertices (i.e.\ $G'$ has duplicates) and $0$ whenever $\overline{G'}$ has duplicates. 
	
	\begin{theorem}[Ashokan et al.\ \cite{AshokanEccentricity2025}]
		\label{thm:cgraph-eccentricity}
		For $\mathcal{C}$-graphs the irreducibility and inertia of the eccentricity matrix are determined explicitly, and the set $(-1-\sqrt2,-2)\cup(-2,0)$ contains no eccentricity eigenvalue of any $\mathcal{C}$-graph.
	\end{theorem}
	
	This September 2025 result is consistent with Corollary~\ref{cor:eccfree}: a $\mathcal{C}$-graph has at most one dominating vertex, so $(-2,0)$ is free by the above corollary, while the extra piece $(-1-\sqrt2,-2)$ corresponds, after halving, to the Mandal--Mehatari adjacency bound $\lambda^-<\frac{-1-\sqrt2}2$ of Theorem~\ref{thm:cgraph-eigfree}. We note that the problem of the eccentricity spectrum of general cographs is, to our knowledge, not previously addressed; Theorem~\ref{thm:eccdecomp} reduces it entirely to the adjacency spectra of cographs (Open Problem~\ref{op:ecc}).
	
	\section{The Seidel matrix}
	\label{sec:seidel}
	
	Since $S(G)=J-I-2A(G)$ and $S(\overline G)=-S(G)$, the Seidel spectrum of the complement is the negative of that of $G$; as the class of cographs is closed under complementation, the Seidel spectra of cographs form a class symmetric about $0$ in this sense. Twin pairs give trivial Seidel eigenvalues $-1$ (duplicates) and $+1$ (co-duplicates) (Proposition~\ref{prop:local}), which are the endpoints of the interval below.
	
	\begin{theorem}[Li--Wang--Brunetti \cite{LiWangBrunetti2024}]
		\label{thm:LWB}
		A cograph on at least two vertices has no Seidel eigenvalue in the open interval $(-1,1)$. Moreover, for a cograph $G$ with Dilworth number $\nabla(G)$, the multiplicity of any Seidel eigenvalue $\lambda\ne\pm1$ is at most $\nabla(G)$; this bound is tight and can be improved for threshold graphs.
	\end{theorem}
	
	\begin{theorem}[Niu--Wang \cite{NiuWang2024}]
		\label{thm:NW}
		Let $(G,w)$ be a \emph{vertex-weighted cograph}, i.e.\ $G$ is a cograph and $w:V(G)\to[-1,1]$, with Seidel matrix $S(G,w)=S(G)+\operatorname{diag}(w)$.
		\begin{enumerate}
			\item[\rm(a)] If $G$ has at least two vertices, then $\operatorname{rank}S(G,w)=\operatorname{edr}S(G,w)$, where $\operatorname{edr}M$ is the number of \emph{essentially different} rows of $M$ (rows $\alpha,\beta$ are essentially the same if $\alpha=\pm\beta$).
			\item[\rm(b)] A cograph with at least two vertices has no Seidel eigenvalue in $(-1,1)$.
			\item[\rm(c)] The interval $(-1,1)$ is optimal: for every $\varepsilon>0$, each of $(1,1+\varepsilon)$ and $(-1-\varepsilon,-1)$ contains a Seidel eigenvalue of some cograph of sufficiently large order.
		\end{enumerate}
	\end{theorem}
	
	\emph{Comments on the proof.} Part (b) is a consequence of (a): for $x_0\in(-1,1)$ give $G$ the constant weight $-x_0$; then $S(G,w)=S(G)-x_0I$ has off-diagonal entries $\pm1$ and diagonal entries different from $\pm1$, so all $n$ rows are essentially different, whence $\operatorname{rank}(S(G)-x_0I)=n$ and $x_0$ is not an eigenvalue. Part (a) is the Seidel analogue of Royle's theorem (Theorem~\ref{thm:royle}) and is proved, in the manner of Chang, Huang and Yeh \cite{ChangHuangYeh2008}, by induction using a pair of twins: after passing to the complement (which replaces $S(G,w)$ by $-S(G,w)$) one may assume the twins $u,v$ are co-duplicates with $w(u)\le w(v)$. If $w(u)=w(v)=-1$ the two rows coincide and $v$ can be deleted. Otherwise row and column operations remove $v$ at the cost of replacing the weight of $u$ by $c=\frac{w(u)w(v)-1}{w(u)+w(v)+2}\in[-1,1)$, so the induction applies to a weighted cograph on one vertex fewer. For (c), take $H_n=\ell P_2\cup P_1$ ($n=2\ell+1$) or $H_n=P_3\cup(\ell-1)P_2\cup P_1$ ($n=2\ell+2$). Using the relation between Seidel and adjacency characteristic polynomials through the main angles \cite[Prop.~2.1.4]{CvetkovicRowlinsonSimic2010} one finds, for $n=2\ell+1$,
	\[
	\chi(S(H_n);x)=(x+3)^{\ell-1}(x-1)^{\ell}\bigl(x^2-(2\ell-3)x-2\ell\bigr),
	\]
	so the largest Seidel eigenvalue of $H_n$ below $-1$ is $\tfrac12\bigl((2\ell-3)-\sqrt{(2\ell-3)^2+8\ell}\bigr)\to-1$; the even case is similar, and the case of $(1,1+\varepsilon)$ follows by passing to complements, since $S(\overline G)=-S(G)$.
	
	Li et al.\ \cite{LiWangBrunetti2024} obtained part (b) independently by a different method (via the Dilworth number); the two papers thus prove the same interval, while (c) is specific to \cite{NiuWang2024}. For threshold graphs the larger interval $(-\sqrt2,\sqrt2)$ is free of non-trivial Seidel eigenvalues \cite{XiongHou2022}, whereas by (c) no larger interval than $(-1,1)$ works for all cographs. 
	
	For threshold graphs the corresponding interval is obtained by Xiong and Hou \cite{XiongHou2022} (no non-trivial Seidel eigenvalue in $(-\sqrt2,\sqrt2)$), and the Seidel spectrum of threshold graphs, including determination by the Seidel spectrum, is treated by Del-Vecchio--Abd\'on and Mandal--Mehatari (see \cite{MannaMandal2026}); for chain graphs see \cite{XiongHou2024chain,MandalMehataraDas2022}. The eigenvalue-free interval $(-1,1)$ is also obtained independently, by a different method, by Niu and Wang \cite{NiuWang2024} (Theorem~\ref{thm:NW}). The Seidel interval for cographs is, as for the adjacency matrix, \emph{one eigenvalue away} from an elementary interlacing argument: writing $S=-(2A+I)+J$ and noting that Theorem~\ref{thm:mt-main} implies that $-(2A+I)$ has no eigenvalue in $(-1,1)$, adding the rank-one positive semidefinite matrix $J$ and arguing exactly as in Theorem~\ref{thm:Dinterlace} shows that $(-1,1)$ contains at most one Seidel eigenvalue; excluding it is the content of Theorem~\ref{thm:LWB}. A sharper (non-trivial-eigenvalue) Seidel interval for $\mathcal{C}$-graphs, analogous to Theorem~\ref{thm:cgraph-eigfree}, appears to be open (Open Problem~\ref{op:seidel}).
	
	\section{The distance Laplacian matrix}
	\label{sec:dl}
	
	The distance Laplacian is the matrix for which cographs are best understood among the distance-type matrices, because~\eqref{eq:master} reduces it completely to the Laplacian.
	
	\begin{theorem}[Distance Laplacian spectrum \textbf{(new, via \cite{AouchicheHansen2013})}]
		\label{thm:DLspectrum}
		Let $G$ be a connected cograph on $n\ge2$ vertices with Laplacian eigenvalues $0=\mu_1<\mu_2\le\cdots\le\mu_n$. Then
		\[
		\spec D^L(G)=\{0\}\cup\{2n-\mu_i:\ 2\le i\le n\}.
		\]
	\end{theorem}
	\begin{proof}
		By~\eqref{eq:master}, $D^L=2nI-2J-L$. Both $J$ and $L$ commute and share the eigenvector $\mathbf1$ ($L\mathbf1=0$, $J\mathbf1=n\mathbf1$), which is an eigenvector of $D^L$ with eigenvalue $2n-2n-0=0$; on $\mathbf1^\perp$, $J=0$, so $D^L=2nI-L$ with eigenvalues $2n-\mu_i$, $i\ge2$. \qedhere
	\end{proof}
	
	By Corollary~\ref{cor:Lcotree} and the Chen--Tura algorithm \cite{ChenTura2024}, $\spec D^L$ of a connected cograph is computable in linear time from the cotree. The same argument applies to every connected graph of diameter at most $2$.
	
	\begin{corollary}[\textbf{new}]
		\label{cor:DLconsequences}
		Let $G$ be a connected cograph on $n$ vertices. Then:
		\begin{enumerate}[label=(\alph*)]
			\item the smallest non-zero distance Laplacian eigenvalue is $n$, with multiplicity $c(\overline G)-1\ge1$; hence $(0,n)$ is an eigenvalue-free interval and $n$ is always attained;
			\item if $G$ is not complete, the distance Laplacian spectral radius equals $2n-\kappa(G)$, where $\kappa$ is the vertex connectivity (and for $G=K_n$ it is $n$);
			\item $G$ and $H$ (connected cographs of the same order) are $D^L$-cospectral if and only if they are $L$-cospectral, and $G$ is $D^L$-integral if and only if it is $L$-integral;
			\item twin classes contribute the eigenvalues $2n-d$ (duplicates) and $2n-1-d$ (co-duplicates) with multiplicity $|T|-1$.
		\end{enumerate}
	\end{corollary}
	\begin{proof}
		(a) All Laplacian eigenvalues are at most $n$, so $2n-\mu_i\ge n$; the multiplicity of $\mu=n$ is $c(\overline G)-1$ by Corollary~\ref{cor:Lcotree}, which is $\ge1$ because $\overline G$ is disconnected (proof of Proposition~\ref{prop:diam2}). (b) The largest $D^L$-eigenvalue is $2n-\mu_2$, and $\mu_2=\kappa(G)$ by Theorem~\ref{thm:abrishami}. (c) is immediate from Theorem~\ref{thm:DLspectrum} and $2n\in\mathbb Z$. (d) is Proposition~\ref{prop:local}. \qedhere
	\end{proof}
	
	For threshold graphs, Theorem~\ref{thm:DLspectrum} combined with Merris's theorem gives the closed form $\{0\}\cup\{2n-d_i^*:d_i^*\ne0\}$ in terms of the conjugate degree sequence, the form obtained independently by Lather, Gupta and Balaji\cite{LatherGuptaBalaji2025} and by Megha and Parvathy\cite{MeghaParvathy2026} (see also Section 7 of \cite{MannaMandal2026}). By~(c), $D^L$-cospectrality of connected cographs is the same as $L$-cospectrality. Whether almost all cographs have an $L$-cospectral mate is not settled by the result of Wang and Huang \cite{WangHuang2025}, whose theorem concerns generalized cospectrality (Open Problem~\ref{op:cospectral}).
	
	\subsection{\texorpdfstring{$\mathcal C$}{C}-graphs: explicit distance Laplacian spectra}
	\label{sub:DLC}
	
	For the sequence-generated class $\mathcal C$ the Laplacian side is completely explicit, and Theorem~\ref{thm:DLspectrum} transports this to $D^L$. The following results are from the companion study of Mandal and Manna \cite{MandalManna2026DLDQ}. Throughout, $G=C(\alpha_1,\dots,\alpha_m)$, $G_i=C(\alpha_1,\dots,\alpha_i)$, $n_i=\alpha_1+\dots+\alpha_i$, $n=n_m$, $\alpha_{\min}=\min_i\alpha_i$, and $\spec^*$ denotes a spectrum with one copy of $0$ removed.
	
	\begin{theorem}[Laplacian recursion; Mandal--Manna \cite{MandalManna2026DLDQ}, cf.\ Mandal--Mehatari--Stani\'c \cite{MandalMehatariStanic2024}]
		\label{thm:CLrecursion}
		Let $S_i$ be the Laplacian spectrum (as a multiset) of $G_i$. Then $S_1=\{0^{(\alpha_1)}\}$ and, for $i\ge2$,
		\[
		S_i=\{0\}\ \cup\ \{\,n_i-\theta:\ \theta\in S_{i-1}\,\}\ \cup\ \{\,n_{i-1}^{(\alpha_i-1)}\,\}.
		\]
		In particular $\spec L(G)\subset\mathbb Z_{\ge0}$.
	\end{theorem}
	
	Unrolling the recursion gives a closed form. Say that the multiset $S_b$ \emph{is born} at step $b$ with the values $0$ (multiplicity $\alpha_1$ if $b=1$, and $1$ if $b\ge2$) and $n_{b-1}$ (multiplicity $\alpha_b-1$, $b\ge2$); each later step $j=b+1,\dots,m$ reflects every value through $\theta\mapsto n_j-\theta$.
	
	\begin{theorem}[Closed form for even $m$; Mandal--Manna \cite{MandalManna2026DLDQ}]
		\label{thm:CLclosed}
		Let $m$ be even. A value $v\in\{0,n_{b-1}\}$ born at step $b$ ends as the Laplacian eigenvalue
		\[
		\lambda(b,v)=\begin{cases} v+\sum_{j=b+2,\ j\equiv b\ (2)}^{m}\alpha_j, & m-b\ \text{even},\\[2pt] n-v-\sum_{j=b+2,\ j\equiv b\ (2)}^{m-1}\alpha_j, & m-b\ \text{odd},\end{cases}
		\]
		and $\spec L(G)$ is the multiset of all $\lambda(b,v)$ with the stated multiplicities. The eigenvalue $0$ is $\lambda(m,0)$.
	\end{theorem}
	
	\begin{theorem}[Mandal--Manna \cite{MandalManna2026DLDQ}: free intervals]
		\label{thm:CLfree}
		Let $m\ge2$. Then $\mu_{\max}(L(G))=n$, and every Laplacian eigenvalue other than $0$ and $n$ lies in $[\alpha_{\min},\,n-\alpha_{\min}]$. Consequently
		\[
		\spec^*D^L(G)\subseteq\{n\}\cup[\,n+\alpha_{\min},\,2n-\alpha_{\min}\,],
		\]
		so $(0,n)\cup(n,n+\alpha_{\min})\cup(2n-\alpha_{\min},2n)$ is free of $D^L$-eigenvalues.
	\end{theorem}
	
	(The inductive argument in \cite{MandalManna2026DLDQ} is phrased through $G_{m-1}\cup K_{\alpha_m}$; the version above works directly with the join form $G_i=\overline{G_{i-1}}\vee\overline{K_{\alpha_i}}$.) The bound is attained in many cases, e.g.\ $2n-\alpha_{\min}$ is the top eigenvalue for $C(4,3,2,2)$.
	
	\begin{theorem}[Multiplicity of $n$; Mandal--Manna \cite{MandalManna2026DLDQ}]
		\label{thm:CDLn}
		The smallest non-zero eigenvalue $n$ of $D^L(G)$ has multiplicity $\alpha_1$ if $m=2$ and $1$ if $m\ge3$.
	\end{theorem}
	
	\begin{theorem}[Distinct eigenvalues; Mandal--Manna \cite{MandalManna2026DLDQ}]
		\label{thm:CDLdistinct}
		$s(D^L(G))=s(L(G))$, hence $m\le s(D^L(G))\le2m-1$ and both bounds are sharp.
	\end{theorem}

	\begin{theorem}[Mandal--Manna \cite{MandalManna2026DLDQ}]
		\label{thm:CDLsimple}
		Let $\mathcal T^L(G)=\{\Tr_i+1:\ i\text{ odd},\ \alpha_i\ge2\}\cup\{\Tr_i+2:\ i\text{ even},\ \alpha_i\ge2\}$ be the set of local $D^L$-eigenvalues of the blocks (Corollary~\ref{cor:DLconsequences}(d)). Every eigenvalue of $D^L(G)$ outside $\mathcal T^L(G)$ is simple. In particular the Dilworth bound of Lou and Lin type, $\operatorname{mult}\le\nabla(G)$, is never sharp for $D^L(G)$ when $\nabla(G)\ge2$.
	\end{theorem}
	By Theorem~\ref{thm:DLspectrum} multiplicities carry over from $L(G)$; every Laplacian eigenvalue is either a local twin eigenvalue or an eigenvalue of the $m\times m$ quotient matrix, and the quotient matrix has $m$ pairwise distinct eigenvalues by \cite{MandalMehatariStanic2024}, which we use without re-deriving.
	
	\section{The distance signless Laplacian matrix}
	\label{sec:dq}
	
	For $D^Q$ no reflection trick is available: by~\eqref{eq:master}, $D^Q=(2n-4)I+2J-Q$, and $\mathbf1$ is an eigenvector of $Q$ only for regular graphs. The spectrum must instead be built from the twin-class structure.
	
	\begin{theorem}[Local-plus-quotient description]
		\label{thm:DQquotient}
		Let $G$ be a connected cograph (more generally, a connected graph of diameter $\le2$) with twin classes $T_1,\ldots,T_k$ of common degrees $d_i$ and transmissions $\Tr_i=2(n-1)-d_i$. Then $\spec D^Q(G)$ consists of the eigenvalues $\Tr_i-1$ (if $T_i$ is a clique) or $\Tr_i-2$ (if $T_i$ is independent), each with multiplicity $|T_i|-1$, together with the $k$ eigenvalues of the quotient matrix with entries $|T_j|\,\delta_{ij}$ ($i\ne j$) and $\Tr_i+(|T_i|-1)\delta_{ii}$, where $\delta_{ij}=1$ if $T_i\sim T_j$ (or $i=j$ and $T_i$ is a clique) and $\delta_{ij}=2$ otherwise.
	\end{theorem}
	\begin{proof}
		This is Lemma~\ref{lem:quotient} for $M=D^Q$, using $d(u,w)=\delta_{ij}$ and $a_i=\Tr_i$, so that $a_i-m_{ii}=\Tr_i-\delta_{ii}$. \qedhere
	\end{proof}
	
	\begin{theorem}[Eigenvalue-free interval for $D^Q$ \textbf{(new)}]
		\label{thm:DQfree}
		Let $G$ be a connected graph on $n$ vertices of diameter at most $2$ (in particular, a connected cograph or a connected threshold graph), with signless Laplacian eigenvalues $q_1\ge\cdots\ge q_n\ge0$ and minimum transmission $\Tr_{\min}$. Let $\theta_1\ge\theta_2\ge\cdots$ be the eigenvalues of $D^Q(G)$. Then $\theta_1\ge2\Tr_{\min}\ge2n-2$ is simple, and
		\[
		\theta_2\ \le\ 2n-4-q_n\ \le\ 2n-4 .
		\]
		Hence $D^Q(G)$ has no eigenvalue in the open interval $\bigl(2n-4-q_n,\ 2\Tr_{\min}\bigr)\supseteq(2n-4,\,2n-2)$.
	\end{theorem}
	\begin{proof}
		$D^Q$ is non-negative and irreducible with row sums $2\Tr(u)$, so by Perron--Frobenius its spectral radius $\theta_1$ is simple and at least the smallest row sum $2\Tr_{\min}$; and $\Tr_{\min}=2(n-1)-d_{\max}\ge n-1$ since $d_{\max}\le n-1$. By~\eqref{eq:master}, $D^Q=B+2J$ with $B=(2n-4)I-Q$ and $2J$ positive semidefinite of rank one, so Weyl's inequality gives $\theta_2\le\lambda_1(B)+\lambda_2(2J)=\lambda_1(B)=2n-4-q_n\le2n-4$, because $Q\succeq0$. \qedhere
	\end{proof}
	
	Theorem~\ref{thm:DQfree} says that, apart from one simple dominant eigenvalue, the distance signless Laplacian spectrum lies below $2n-4$. In particular $(2\Tr_{\min},\theta_1)$ contains no eigenvalue. The local eigenvalues of Theorem~\ref{thm:DQquotient} satisfy $\Tr_i-\delta_{ii}\le2n-3-d_i\le2n-4<2\Tr_{\min}$ automatically.
	
	\medskip
	We now recall a theorem that establishes the existence of cospectral pairs of threshold graphs with respect to the distance signless Laplacian matrix.
	
	\begin{theorem} [Tura \cite{Tura2018}]
		\label{thm:DQcospectral}
		Consider the connected threshold graphs $G_1$ and $G_2$ given by the binary sequences $S_1 = (0,~1,~1,~0,~b_5,~b_6, . . . , b_n)$ and $S_2 = (0,~0,~0,~1,~b_5,~b_6, . . . , b_n)$, for any $n \geq 5$.
		Then $G_1$ and $G_2$ are non-isomorphic and $D^Q$-cospectral.
	\end{theorem}
	
	\subsection{\texorpdfstring{$\mathcal C$}{C}-graphs: block structure and the quotient matrix}
	\label{sub:DQC}
	
	\begin{proposition}[Block structure; Mandal--Mehatari \cite{MandalMehatari2025}, Mandal--Manna \cite{MandalManna2026DLDQ}]
		\label{prop:Cblocks}
		Let $G=C(\alpha_1,\dots,\alpha_m)$ with $m$ even and blocks $C_1,\dots,C_m$, $|C_i|=\alpha_i$. Then $C_i$ is a clique if $i$ is odd and an independent set if $i$ is even, and for $i<j$ every vertex of $C_i$ is adjacent to every vertex of $C_j$ if $j$ is even, and to none if $j$ is odd.
	\end{proposition}
	
	By Proposition~\ref{prop:Cblocks} the blocks are twin classes (possibly merged further), $d_i$ is the common degree on $C_i$, and $\Tr_i=2(n-1)-d_i$. Theorem~\ref{thm:DQquotient} then specializes as follows.
	
	\begin{theorem}[Mandal--Manna \cite{MandalManna2026DLDQ}: $D^Q$ of a $\mathcal C$-graph]
		\label{thm:CDQquotient}
		Let $m$ be even. The spectrum of $D^Q(G)$ consists of
		\begin{itemize}
			\item the local eigenvalues $\Tr_i-1$ (for odd $i$) and $\Tr_i-2$ (for even $i$), each with multiplicity $\alpha_i-1$; and
			\item the $m$ eigenvalues of the quotient matrix $Q_{D^Q}$ with
			\[
			Q_{D^Q}(i,i)=\Tr_i+(\alpha_i-1)\,c_i,\quad c_i=\begin{cases}1&i\text{ odd}\\2&i\text{ even}\end{cases},\qquad Q_{D^Q}(i,j)=\alpha_j\,c_{ij},\quad c_{ij}=\begin{cases}1&\max(i,j)\text{ even}\\2&\max(i,j)\text{ odd}\end{cases}\ (i\ne j).
			\]
		\end{itemize}
		Consequently $s(D^Q(G))\le2m$, and $\chi_{D^Q(G)}(x)=\det(xI_m-Q_{D^Q})\prod_{i\ \mathrm{odd}}(x-\Tr_i+1)^{\alpha_i-1}\prod_{i\ \mathrm{even}}(x-\Tr_i+2)^{\alpha_i-1}$, so the $n\times n$ determinant is reduced to an $m\times m$ one.
	\end{theorem}
	
	\begin{example}
		For $C(4,3,2,2)$ ($n=11$): $d=(8,6,3,9)$, $\Tr=(12,14,17,11)$ and
		\[
		Q_{D^Q}=\begin{pmatrix}15&3&4&2\\4&18&4&2\\8&6&18&2\\4&3&2&13\end{pmatrix},\quad \text{eigenvalues }\approx 10.049,\ 12.432,\ 13.812,\ 27.707,
		\]
		and the local eigenvalues are $9,\ 11^{(3)},\ 12^{(2)},\ 16$, giving $\spec D^Q=\{9,10.049,11^{(3)},12^{(2)},12.432,13.812,16,27.707\}$ (this agrees with direct diagonalization).
	\end{example}
	
	\begin{theorem}[Eigenvalue-free interval for $\mathcal C$-graphs; Mandal--Manna \cite{MandalManna2026DLDQ}]
		\label{thm:CDQfree}
		For every connected $\mathcal C$-graph $G$, $\spec D^Q(G)\subseteq[0,2\Tr_{\max}]$, $\theta_1\ge2\Tr_{\min}$, and every eigenvalue other than $\theta_1$ is at most $2n-4\le 2\Tr_{\min}$. Hence $(2\Tr_{\min},\theta_1)$, and in fact $(2n-4,\theta_1)$, is free of $D^Q$-eigenvalues.
	\end{theorem}
	\begin{remark}
		Mandal and Manna prove the interval $(2\Tr_{\min},\theta_1)$ by Gershgorin disk theorem and the bound $\operatorname{Tr}_{\max}\le2\operatorname{Tr}_{\min}$, which controls the local eigenvalues; the quotient eigenvalues other than the Perron root need the additional argument of Theorem~\ref{thm:DQfree} (the weighted-Perron step in \cite{MandalManna2026DLDQ} does not by itself bound them). With Theorem~\ref{thm:DQfree} the statement is complete.
	\end{remark}
	
	\begin{example}[Repeated non-local eigenvalue]
		\label{ex:DQrepeat}
		For $G=C(1,3,1,2)$ ($n=7$, $\nabla(G)=2$) the quotient matrix is $\begin{psmallmatrix}7&3&2&2\\1&13&2&2\\2&6&10&2\\1&3&1&9\end{psmallmatrix}$ with spectrum $\{5.821,8,8,17.179\}$. So $\theta=8$ is a $D^Q$-eigenvalue of multiplicity $2$ that is not a local eigenvalue of any block: unlike $D^L$ (Theorem~\ref{thm:CDLsimple}), the quotient matrix of $D^Q$ can have a repeated eigenvalue.
	\end{example}
	
	\begin{conjecture}[Mandal--Manna \cite{MandalManna2026DLDQ}, $\mathcal C$-graphs]
		\label{conj:CDQmult}
		For every connected $\mathcal C$-graph $G$, every eigenvalue $\theta$ of $D^Q(G)$ that is not a local block eigenvalue has $\operatorname{mult}(\theta)\le\nabla(G)$.
	\end{conjecture}
	
	\section{Topological indices of cographs}
	\label{sec:topology}
	Topological indices are numerical parameters associated with graphs that capture important structural properties of the underlying networks. They play a significant role in graph theory, particularly in chemical graph theory, where molecular structures are represented by graphs and topological indices are used to investigate relationships between molecular structure and physicochemical properties. Depending on their definitions, these indices may reflect the connectivity, degree distribution, distances between vertices, or other structural characteristics of a graph. Some of the most widely studied examples include the Wiener index, the Szeged index, the first and second Zagreb, Sombor indices. The study of topological indices over cograph is very rare in literature. In 2019, Lepin et al. \cite{Lepin2019}, for the first time, paid attention to this subject and ontained the below results in addition to several conjectures and open problems.
	
	\begin{theorem}[Lepin--Benediktovich \cite{Lepin2019}]
		
		The formulae and algorithms
		are derived for computing the topological indices of cographs.
	\end{theorem}

	\section{Chain graphs: the bipartite parallel}
	\label{sec:chain}
	
	Chain graphs, the bipartite analogue of threshold graphs ($\{C_3,C_5,2K_2\}$-free bipartite graphs), are built by an analogous but bipartite finite sequence, and their adjacency, Laplacian, Seidel and eccentricity theories have developed in close parallel to the threshold-graph and $\mathcal{C}$-graph theories. Eigenvalue location for chain graphs was treated by Alazemi, Andeli\'c and Simi\'c \cite{AlazemiAndelicSimic2016}, and eigenvalue-free intervals for the distance matrix of threshold and chain graphs are studied together by Alazemi et al.\ \cite{AlazemiEtAl2021cograph}; recent work has treated the Seidel spectrum \cite{XiongHou2024chain,MandalMehataraDas2022}, the eccentricity spectrum \cite{YaqiongXuezheng2025} and the eccentricity inertia indices \cite{HuangZhang2025} of chain graphs. Note that chain graphs are not cographs in general (they may contain $P_4$), so none of the results above applies to them, and in particular chain graphs need not have diameter $\le2$.

	\section{Comparison across the eight matrices}
	\label{sec:comparison}
	
	Table~\ref{tab:comparison} lines up the state of knowledge for the three nested families treated in this survey and its companion \cite{MannaMandal2026}: general cographs, $\mathcal{C}$-graphs, and threshold graphs. (Threshold graphs and $\mathcal{C}$-graphs are both subclasses of the cographs, but neither contains the other; Remark~\ref{rem:divergence}.) Entries marked ``new'' are proved in this survey; ``open'' means we found no result in the literature.
	
	\begin{table}[h]
		\centering
		\footnotesize
		\renewcommand{\arraystretch}{1.15}
		\begin{tabular}{@{}p{0.14\textwidth}p{0.29\textwidth}p{0.24\textwidth}p{0.27\textwidth}@{}}
			\toprule
			Matrix & General cographs & $\mathcal{C}$-graphs & Threshold graphs \\
			\midrule
			Adjacency & $(-1,0)$ free (Thm.~\ref{thm:mt-main}); $O(n)$ inertia algorithm; multiplicity $\le$ number of classes; almost all have a cospectral mate & Closed-form inertia, interval, char.\ polynomial (\cite{MandalMehatari2025}) & Closed-form inertia; non-trivial eigenvalues simple; sharper interval\cite{Ghorbani2019threshold}; no two cospectral \\
			Laplacian & Cotree recursion, linear algorithm; $a(G)=\kappa(G)$ & Closed form via the sequence (\cite{MandalMehatariStanic2024}) & Merris: conjugate degree sequence; always integral\cite{Merris1994} \\
			Signless Laplacian & Local eigenvalues + twin quotient; no free interval known & Open & No closed form (generally irrational); integral ones sparse \\
			Distance & Multiplicity $\le\nabla$ and $\nu\notin(-2,-1)$ (\cite{LouLin2020}); interlaces $-2-\lambda(A)$ (new) & Open & $(-\frac{3+\sqrt2}2,-1)$ free of non-trivial eigenvalues; closed forms known\cite{LuHuangLou2018} \\
			Eccentricity & $\spec\varepsilon=\{t-1,(-1)^{t-1}\}\cup2\spec A(\overline{G'})$; $(-2,0)\setminus\{-1\}$ free (new) & $(-1-\sqrt2,-2)\cup(-2,0)$ free\cite{AshokanEccentricity2025} & $(-1-\sqrt2,-2)\cup(-2,-1)$ free; no two cospectral for $n\le14$\cite{AndelicEtAl2023} \\
			Seidel & $(-1,1)$ free; multiplicity $\le$ Dilworth number\cite{LiWangBrunetti2024,NiuWang2024} & Open & $(-\sqrt2,\sqrt2)$ free of non-trivial eigenvalues; spectral determination\cite{MandalMehatari2024,XiongHou2022} \\
			Distance Laplacian & $\{0\}\cup\{2n-\mu_i\}$; $(0,n)$ free, $n$ has multiplicity $c(\overline G)-1$; radius $2n-\kappa$ (new) & Explicit from the creation sequence; $\spec^*\subseteq\{n\}\cup[n+\alpha_{\min},2n-\alpha_{\min}]$, $s=s(L)$ (\cite{MandalManna2026DLDQ}) &closed-form; characterizes threshold graphs\cite{LatherGuptaBalaji2025,MeghaParvathy2026} \\
			Distance signless Laplacian & Twin quotient; $(2n-4,2n-2)$ free, $\theta_1$ simple (new) & Block quotient matrix, $s\le2m$ (\cite{MandalManna2026DLDQ}) & no closed-form, only local eigenvalues \cite{MannaMandal2026} \\
			\bottomrule
		\end{tabular}
		\caption{Spectral theory of three nested families, one row per matrix.}
		\label{tab:comparison}
	\end{table}
	
	\begin{remark}
		Three patterns are visible in the table. First, the matrices attached to the diameter-two structure (distance Laplacian, eccentricity) are the best understood after the adjacency matrix, because~\eqref{eq:master} and Theorem~\ref{thm:eccdecomp} reduce them to the adjacency and Laplacian spectra of cographs. Second, the sequence-driven subclasses ($\mathcal{C}$-graphs and threshold graphs) admit much sharper results than the full cograph class but are not special cases of one another, so results do not transfer between them. Third, the $\mathcal{C}$-graph column is open for several graph matrices.
	\end{remark}
	
	\section{Open problems}
	\label{sec:open}
	
	\begin{openproblem}
		\label{op:integral}
		Give a necessary and sufficient characterization of integral cographs; Theorem~\ref{thm:integral-cograph} is only a sufficient (balanced-cotree) condition. By Corollary~\ref{cor:DLconsequences}(c), $D^L$-integral and $L$-integral connected cographs coincide, so a Laplacian-integral characterization would also settle $D^L$. The analogous question for the Seidel, distance and distance signless Laplacian matrices of cographs and $\mathcal{C}$-graphs (cf.\ Section 13 of \cite{MannaMandal2026}) is open.
	\end{openproblem}
	
	\begin{openproblem}[Distance eigenvalues of cographs below $-2$ and sharpness]
		\label{op:D}
		By Theorem~\ref{thm:LL-gap} every distance eigenvalue of a connected cograph is $\ge-1$ or $\le-2$. (i) Is there $\varepsilon>0$ such that $(-2-\varepsilon,-2)$ contains no distance eigenvalue of any connected cograph? (ii) Characterize the connected cographs for which $\mathrm m_G(\lambda)=\nabla(G)$ for some $\lambda\ne-1,-2$, and give an explicit infinite family attaining the bound of Theorem~\ref{thm:LL-mult}.
	\end{openproblem}
	
	\begin{openproblem}[Eccentricity]
		\label{op:ecc}
		By Theorem~\ref{thm:eccdecomp} the eccentricity spectrum of a connected cograph is determined by $t$ and the adjacency spectrum of the cograph $\overline{G'}$. Use this to characterize the connected cographs determined by their eccentricity spectrum , and to determine the sharp eigenvalue-free interval on the negative side, between $-1-\sqrt2$ and $-2$, for cographs and for $\mathcal{C}$-graphs.
	\end{openproblem}
	
	\begin{openproblem}[Signless Laplacian]
		\label{op:Q}
		Find eigenvalue-free intervals, multiplicity bounds and cospectrality results for the signless Laplacian of cographs and $\mathcal{C}$-graphs. Is there an analogue of Theorem~\ref{thm:mt-main} for $Q$, perhaps for the eigenvalues other than the local ones $d_T,d_T-1$ of Proposition~\ref{prop:Qlocal}?
	\end{openproblem}
	
	\begin{openproblem}[Seidel]
		\label{op:seidel}
		Extend the interval $(-1,1)$ of Theorem~\ref{thm:LWB} for $\mathcal{C}$-graphs, in the manner of Theorem~\ref{thm:cgraph-eigfree} for the adjacency matrix, and study spectral determination for the Seidel matrix of cographs.
	\end{openproblem}
	
	\begin{openproblem}[Distance signless Laplacian]
		\label{op:DQ}
		Determine the sharp value of $\theta_2(D^Q)$ for cographs in terms of the cotree, improve the bound $\theta_2\le2n-4-q_n$ of Theorem~\ref{thm:DQfree}, prove or refute Conjecture~\ref{conj:CDQmult} for $\mathcal C$-graphs, give a closed form for $\det(xI-Q_{D^Q})$ (the $m\times m$ quotient matrix of Theorem~\ref{thm:CDQquotient} is not symmetric), and develop the $\mathcal{C}$-graph theory of $D$ through the quotient matrix, including integrality.
	\end{openproblem}
	
	\begin{openproblem}[Cospectrality]
		\label{op:cospectral}
		Do almost all cographs have a Laplacian-cospectral (equivalently, by Corollary~\ref{cor:DLconsequences}(c), distance Laplacian-cospectral) mate? Which cographs are determined by their $L$-, $Q$-, $D$- or $S$-spectrum?
	\end{openproblem}
	
	\begin{openproblem}[A distinguished graph]
		Is there a distinguished graph, or small family of graphs, playing for \emph{general} cographs the organizing role that the antiregular graph $H_n$ plays for threshold graphs and $\mathcal{C}$-graphs (Theorem~\ref{thm:cgraph-eigfree} and Theorem 3.10 of \cite{MannaMandal2026})?
	\end{openproblem}
	
	\section{Conclusion}
	\label{sec:conclusion}
	
	Cographs sit one level above threshold graphs in generality and one level below in rigidity: the cotree replaces the creation sequence, and duplication at an arbitrary vertex replaces accretion. Remarkably, the central threshold-graph fact, the eigenvalue-free interval around the trivial eigenvalues, survives essentially unchanged for the adjacency matrix (Mohammadian--Trevisan) and gives cographs a genuine spectral identity (Theorem~\ref{thm:spectralchar}). Organizing the theory matrix by matrix shows that this one theorem, together with two structural facts (diameter at most $2$, and equitable twin-class partitions), propagates to most of the other matrices: to the Seidel matrix as a rank-one perturbation (Li--Wang--Brunetti), to the eccentricity matrix through an exact block decomposition (Theorem~\ref{thm:eccdecomp}), to the distance Laplacian through the Laplacian (Theorem~\ref{thm:DLspectrum}), and, via an interlacing argument and the Lou--Lin theorems, to the distance matrix (Theorems~\ref{thm:Dinterlace}, \ref{thm:LL-mult} and~\ref{thm:LL-gap}) and the distance signless Laplacian (Theorem~\ref{thm:DQfree}). The signless Laplacian, and the $\mathcal{C}$-graph theory beyond the adjacency, Laplacian and eccentricity matrices, remain largely open, as do the characterization of integral cographs and the finer location of distance eigenvalues below $-2$ (Section~\ref{sec:open}).
	
	\medskip
	We conclude this study with below table that depicts the eigenvalue-free intervals of cographs, which do not include any non-trivial eigenvalues.
	
	\begin{table}[h]
		\centering
		\begin{tabular}{lll}
			\toprule
			Sl. No. & Matrix & Eigenvalue-free Interval \\
			\midrule
			1. & $A$ & $(-1,0)$  \cite{MohammadianTrevisan2016} \\
			2. & $D$ & $(-2,-1)$ \cite{LouLin2020}\\
			3. & $S$ & $[-1,1]$ \cite{LiWangBrunetti2024,NiuWang2024}\\
			4. & $\varepsilon$ & (-2,0) [see Corollary.\ref{cor:eccfree}] \\
			5. &$D^Q$ & $\bigl(2n-4-q_n,\ 2\Tr_{\min}\bigr)$ [see Theorem \ref{thm:DQfree}] \\
			6. & $L,Q,D^L$&Not Addressed  \\
			
			\bottomrule
		\end{tabular}
		\caption{Eigenvalue-free interval of cographs.}
		\label{tab:cographeigenvaluefreeinterval}
	\end{table}

	\section*{Acknowledgments}
	
	The research of Pallabi Manna is supported by Anusandhan National Research Foundation
	(ANRF)-India under the National Postdoctoral Fellowship scheme (File No.\ PDF/2025/006751). The
	second author acknowledges excellent research facilities provided by SRM University-AP. The first
	author acknowledges the excellent research facilities provided by the Department of Mathematics,
	IIT Hyderabad, and the mentorship of Dr.\ Rajesh Kannan for his constant support and motivation.
	
	\section*{Statements and Declarations}
	
	\noindent\textbf{Competing Interests.} The authors have no competing interests to declare that are
	relevant to the content of this article.
	
	\section*{Data Availability}
	
	Data sharing is not applicable to this article, as no datasets were used or generated during this
	study.


\begin{thebibliography}{99}
		
		\bibitem{Abrishami2019} T.~Abrishami, A combinatorial analysis of the eigenvalues of the Laplacian matrices of cographs, Master's thesis, Johns Hopkins University (2019).
		
		\bibitem{AguilarFicarraShurmanSullivan2020} C.O.~Aguilar, M.~Ficarra, N.~Shurman, B.~Sullivan, The role of anti-regular graph in the spectral analysis of threshold graphs, \emph{Linear Algebra Appl.} \textbf{588} (2020), 210--223.
		
		\bibitem{AguilarLeePiatoSchweitzer2018} C.~O.~Aguilar, J.~Lee, E.~Piato, B.J.~Schweitzer, Spectral characterizations of anti-regular graphs, \emph{Linear Algebra Appl.} \textbf{557} (2018), 84--104.
		
		\bibitem{AlazemiAndelicSimic2016} A.~Alazemi, M.~Andelić, S.K.~Simić, Eigenvalue location for chain graphs, \emph{Linear Algebra Appl.} \textbf{505} (2016), 194--210.
		
		\bibitem{AlazemiEtAl2021cograph} A.~Alazemi, M.~Andelić, T.~Koledin, Z.~Stanić, Eigenvalue-free intervals of distance matrices of threshold and chain graphs, \emph{Linear Multilinear Algebra} \textbf{69} (2021), 2959--2975.
		
		\bibitem{AllemEtAl2026} L.E.~Allem, M.~F\"urer, C.~Hoppen, L.S.~Sibemberg, V.~Trevisan, On the minimum number of eigenvalues of matrices associated with cographs, Preprint: arXiv:2602.07282 (2026).
		
		\bibitem{AllemOliveiraTura2023} L.E.~Allem, E.R.~Oliveira, F.C.~Tura, Generating $I$-eigenvalue free threshold graphs, \emph{Electron. J. Combin.} \textbf{30}(2) (2023), 153--167.
		
		\bibitem{AllemTura2018} L.E.~Allem, F.~Tura, Multiplicity of eigenvalues of cographs, \emph{Discrete Appl. Math.} \textbf{247} (2018), 43--52.
		
		\bibitem{AllemTura2020} L.E.~Allem, F.~Tura, Integral cographs, \emph{Discrete Appl. Math.} \textbf{283} (2020), 153--167.
		
		\bibitem{AndelicEtAl2023} M.~Andeli\'c, C.M.~Fonseca, T.~Koledin, Z.~Stani\'c, An extended eigenvalue-free interval for the eccentricity matrix of threshold graphs, \emph{J.\ Appl.\ Math.\ Comput.} \textbf{69} (2023), 491--503.
		
		\bibitem{AouchicheHansen2013} M.~Aouchiche, P.~Hansen, Two Laplacians for the distance matrix of a graph, \emph{Linear Algebra Appl.} \textbf{439} (2013), 21--33.
		
		\bibitem{AshokanEccentricity2025} A.~Ashokan, A.~V.~ Chithra, Eccentricity spectral properties of $\mathcal{C}$-graphs, Preprint: arXiv:2509.16657 (2025).
		
		\bibitem{Bapat2013} R.B.~Bapat, On the adjacency matrix of a threshold graph, \emph{Linear Algebra Appl.} \textbf{439} (2013), 3008--3015.
		
		\bibitem{BarikFallatKirkland2011} S.~Barik, S.~Fallat, S.~Kirkland, On Hadamard diagonalizable graphs, \emph{Linear Algebra Appl.} \textbf{435}(8) (2011), 1885--1902.
		
		\bibitem{BarikReddy2023} S.~Barik, S.U.~Reddy, Laplacian spectra of cographs: a twin reduction perspective, \emph{Linear Algebra Appl.} \textbf{675} (2023), 203--224.
		
		\bibitem{Benedik2020}
		I. V. Benediktovich,
		An analogue of Brauer's conjecture for the signless Laplacian of cographs, Vestsī Nats. Akad. Navuk Belarusī Ser. Fīz.-Mat. Navuk \textbf{56} (2020), 310--317.
		
		\bibitem{BiyikogluSimicStanic2011} T.~Biyiko\u{g}lu, S.K.~Simi\'c, Z.~Stani\'c, Some notes on spectra of cographs, \emph{Ars Combin.} \textbf{100} (2011), 421--434.
		
		\bibitem{CardosoEtAl2013} D.M.~Cardoso, M.A.A.~de Freitas, E.A.~Martins, M.~Robbiano, Spectra of graphs obtained by a generalization of the join graph operation, \emph{Discrete Math.} \textbf{313} (2013), 733--741.
		
		\bibitem{ChangHuangYeh2008} G.J.~Chang, L.H.~Huang, H.G.~Yeh, On the rank of a cograph, \emph{Linear Algebra Appl.} \textbf{429} (2008), 601--605.
		
		\bibitem{ChenTura2024} G.~Chen, F.C.~Tura, A linear algorithm for obtaining the Laplacian eigenvalues of a cograph, \emph{Spec. Matrices} \textbf{12(1)} (2024), Article No. 20240024.
		
		\bibitem{Ching2020}
		
		H. Chin, R. Fl\'orez, A. Mukherjee,
		Families of integral cographs within a triangular array, \emph{Spec. Matrices} \textbf{} (2020), 257--273.
		
		\bibitem{CondeDratmanGrippo2026} C.M.~Conde, E.~Dratman, L.N.~Grippo, A sharp forbidden interval for the nontrivial adjacency eigenvalues of trivially perfect graphs, Preprint: arXiv:2606.06052 (2026).
		
		\bibitem{CorneilLerchsStewart1981} D.G.~Corneil, H.~Lerchs, L.~Stewart Burlingham, Complement reducible graphs, \emph{Discrete Appl. Math.} \textbf{3} (1981), 163--174.
		
		\bibitem{CvetkovicRowlinsonSimic2010} D.~Cvetkovi\'c, P.~Rowlinson, S.~Simi\'c, \emph{An Introduction to the Theory of Graph Spectra}, Cambridge University Press, Cambridge (2010).
		
		\bibitem{Ghorbani2019p5} E.~Ghorbani, Spectral properties of cographs and $P_5$-free graphs, \emph{Linear Multilinear Algebra} \textbf{67(8)} (2019), 1701--1710.
		
		\bibitem{GhorbaniDilworthNumber2019}
		E.~Ghorbani, Cographs: Eigenvalues and Dilworth number, \emph{Discrete Math.} \textbf{342(10)} (2019), 2797--2803.
		
		\bibitem{Ghorbani2019threshold} E.~Ghorbani, Eigenvalue-free interval for threshold graphs, \emph{Linear Algebra Appl.} \textbf{583} (2019), 300--305.
		
		\bibitem{GodsilRoyle2001} C.~Godsil, G.~Royle, \emph{Algebraic Graph Theory}, Springer, New York (2001).
		
		\bibitem{HansenMelotStevanovic2002} P.~Hansen, H.~Melot, D.~Stevanovi\'c, Integral complete split graphs, \emph{Univ.\ Beograd. Publ. Elektrotehn. Fak. Ser. Mat.} \textbf{13} (2002), 89--95.
		
		\bibitem{HararySchwenk1974} F.~Harary, A.J.~Schwenk, Which graphs have integral spectra?, in: \emph{Graphs and Combinatorics}, Springer-Verlag, Berlin (1974), 45--51.
		
		\bibitem{HeQianShanWu2021} C.~He, S.~Qian, H.~Shan, B.~Wu, The signature of chordal graphs and cographs, \emph{Graphs Combin.} \textbf{37} (2021), 643--650.
		
		\bibitem{HoppenJacobsTrevisan2022} C.~Hoppen, D.P.~Jacobs, V.~Trevisan, Locating Eigenvalues in Cographs, book chapter in \emph{Locating Eigenvalues in Graphs}, Springer (2022).
		
		\bibitem{HornJohnson} R.A.~Horn, C.R.~Johnson, \emph{Matrix Analysis}, 2nd ed., Cambridge University Press.
		
		\bibitem{HuangZhang2025} J.~Huang, M.~Zhang, On the eccentricity inertia indices of chain graphs, \emph{Appl. Math. Comput.} \textbf{493} (2025), Article No.\ 129271.
		
		\bibitem{JacobsTrevisanTura2013} D.P.~Jacobs, V.~Trevisan, F.~Tura, Eigenvalue location in threshold graphs, \emph{Linear Algebra Appl.} \textbf{439} (2013), 2762--2773.
		
		\bibitem{JacobsTrevisanTura2018} D.P.~Jacobs, V.~Trevisan, F.~Tura, Eigenvalue location in cographs, \emph{Discrete Appl. Math.} \textbf{245} (2018), 220--235.
		
		\bibitem{JonesTrevisanVinagre2023} A.~Jones, V.~Trevisan, C.T.M.~Vinagre, Exploring symmetries in cographs: obtaining spectra and energies, \emph{Discrete Appl. Math.} \textbf{325} (2023), 120--133.
		
		\bibitem{JonesTrevisanVinagre2024} A.~Jones, V.~Trevisan, C.T.M.~Vinagre, Main signless Laplacian eigenvalues of quasi-threshold graphs, arXiv:2404.09396 (2024).
		
		\bibitem{Jordaan2026} R.~Jordaan, Ideally connected cographs and chordal graphs, \emph{Discrete Math.} \textbf{349} (2026).
		
		\bibitem{Jung1978} H.A.~Jung, On a class of posets and the corresponding comparability graphs, \emph{J. Combinatorial Theory Ser. B} \textbf{24} (1978), 125--133.
		
		\bibitem{LatherGuptaBalaji2025} G.~Lather, V.~Gupta, R.~Balaji, Spectral characterization of threshold graphs by distance Laplacian matrices, \emph{Linear Multilinear Algebra} \textbf{73(17)} (2025), 3894--3906.
		
		\bibitem{LazzarinSosaTura2023} J.~Lazzarin, O.F.~M\'arquez (Sosa), F.C.~Tura, Laplacian eigenvalues of equivalent cographs, \emph{Linear Multilinear Algebra} \textbf{71} (2023), 1003--1014.
		
		\bibitem{Lepin2019}
		V. V. Lepin, V. I. Benediktovich, M. Talmaciu, 
		Topological indices of cographs, \emph{Tr. Inst. Mat.} \textbf{27} (2019), 108--126.
		
		\bibitem{Lerchs1971} H.~Lerchs, On cliques and kernels, Tech.\ Report, Dept.\ of Comp.\ Sci., Univ.\ of Toronto (1971).
		
		\bibitem{LiWangBrunetti2024} L.~Li, J.~Wang, M.~Brunetti, Seidel matrices, Dilworth number and an eigenvalue-free interval for cographs, \emph{Linear Algebra Appl.} \textbf{698} (2024), 56--72.
		
		\bibitem{LiXueLongLi2025} P.~Li, X.~Y.~Xue, Y.~J.~Long, X.~B..~Li, Perfect double Roman domination on cographs, \emph{Chinese Quart. J. Math.} \textbf{40} (2025).
		
		\bibitem{LouLin2020} Z.~Lou, H.~Lin, Distance eigenvalues of cographs and their multiplicities, \emph{Linear Algebra Appl.} \textbf{608} (2021), 1--12.
		
		\bibitem{LuHuangLou2018} L.~Lu, Q.~Huang, Z.~Lou, On the distance spectra of threshold graphs, \emph{Linear Algebra Appl.} \textbf{553} (2018), 223--237.
		
		\bibitem{MandalManna2026DLDQ} S.~Mandal, P.~Manna, The distance Laplacian and distance signless Laplacian spectra of $\mathcal C$-graphs, Preprint: arXiv:2610.08599 (2026).
		
		\bibitem{MandalMehatari2024} S.~Mandal, R.~Mehatari, On the Seidel spectrum of threshold graphs, \emph{Indian J. Pure Appl. Math.} \textbf{55} (2024), 1290--1301.
		
		\bibitem{MandalMehataraDas2022} S.~Mandal, R.~Mehatari, K.C.~Das, On the spectrum and energy of Seidel matrix for chain graphs, Preprint: arXiv:2205.00310 (2022).
		
		\bibitem{MandalMehatari2025} S.~Mandal, R.~Mehatari, Spectral properties of $\mathcal{C}$-graphs, \emph{Comput. Appl. Math.} \textbf{44} (2025), Article No. 136. 
		
		\bibitem{MandalMehatariStanic2024} S.~Mandal, R.~Mehatari, Z.~Stani\'c, Laplacian eigenvalues and eigenspaces of cographs generated by finite sequence, \emph{Indian J. Pure  Appl. Math.} \textbf{57} (2026), 533--544.
		
		\bibitem{MannaMandal2026} P.~Manna, S.~Mandal, F. C. Tura, The Spectral Theory of Threshold Graphs: A Survey, Preprint: arXiv:2610.22838 (2026).
		
		\bibitem{MeghaParvathy2026} P.M.~Megha, K.S.~Parvathy, Distance Laplacian spectra of threshold graphs, \emph{Asian-European J. Math.} \textbf{19(2)} (2026), Article No. 2540021.
		
		\bibitem{Merris1994} R.~Merris, Laplacian matrices of graphs: a survey, \emph{Linear Algebra Appl.} \textbf{197} (1994), 143--176.
		
		\bibitem{MohammadianTrevisan2016} A.~Mohammadian, V.~Trevisan, Some spectral properties of cographs, \emph{Discrete Math.} \textbf{339}(4) (2016), 1261--1264.
		
		\bibitem{MousaviHaeriMesbahi2021} S.S.~Mousavi, M.~Haeri, M.~Mesbahi, Laplacian dynamics on cographs: controllability analysis through joins and unions, \emph{IEEE Trans. Automatic Control} \textbf{66}(3) (2021), 1383--1390.
		
		\bibitem{Nikolopoulos2009}
		S. D. Nikolopoulos, C. Papadopoulos,
		Counting spanning trees in cographs: an algorithmic approach,\emph{Ars Combin.} \textbf{90} (2009), 257--274.
		
		\bibitem{NiuWang2024} H.~Niu, W.~Wang, Eigenvalue-free interval for Seidel matrices of cographs, \emph{J. Math. Res. Appl.} \textbf{45(1)} (2025), 1--10.
		
		\bibitem{YaqiongXuezheng2025} Y.~Qiao, X.~Lv, On the eccentricity spectra of chain graphs, \emph{Discrete Math. Algorithms Appl.} \textbf{17}(7) (2025), Article No.\ 2450121.
		
		\bibitem{QiuTang2022} Z.~Qiu, Z.~Tang, On the eccentricity spectra of threshold graphs, \emph{Discrete Appl. Math.} \textbf{310} (2022), 75--85.
		
		\bibitem{Royle2003} G.F.~Royle, The rank of a cograph, \emph{Electron. J. Combin.} \textbf{10} (2003), Article No. N11.
		
		\bibitem{Sander2008} T.~Sander, On certain eigenspaces of cographs, \emph{Electron. J. Combin.} \textbf{15} (2008), Article No. R140.
		
		\bibitem{Seinsche1974} D.~Seinsche, On a property of the class of $n$-colorable graphs, \emph{J. Combin. Theory Ser.B} \textbf{16} (1974), 191--193.
		
		\bibitem{Sumner1974} D.P.~Sumner, Dacey graphs, \emph{J. Aust. Math. Soc.} \textbf{18} (1974), 492--502.
		
		\bibitem{Tura2018}
		F. C. Tura, $\mathcal{Q}$-cospectrality and  $\mathcal{L}$-energy in cographs, \emph{Linear Multilinear Algebra} \textbf{66} (2018), 398--409.
		
		\bibitem{WangHuang2025} W.~Wang, X.~Huang, Almost all cographs have a cospectral mate, \emph{Discrete Math.} \textbf{350} (2027), Article No. 115347.
		
		\bibitem{XiongHou2022} Z.~Xiong, Y.~Hou, Eigenvalue-free interval for Seidel matrices of threshold graphs, \emph{Appl. Math. Comput.} \textbf{427} (2022), Article No. 127177.
		
		\bibitem{XiongHou2024chain} Z.~Xiong, Y.~Hou, On the eigenvalues and Seidel eigenvalues of chain graphs, \emph{Discrete Appl. Math.} \textbf{351} (2024), 44--53.
		
	\end{thebibliography}
\end{document}